\documentclass[11pt]{article}
\input amssym.def
\usepackage{graphicx}
\usepackage[latin1]{inputenc}

\def\sqw{\hbox{\rlap{\leavevmode\raise.3ex\hbox{$\sqcap$}}$%
\sqcup$}}

  \def\bn{\hbox{\it 
I\hskip -2pt N}}      
\def\demo{\noindent{\bf Proof \ }} \newtheorem{theorem}{Theorem} 
\newtheorem{lemma}{Lemma} \newtheorem{proposition}{Proposition} 
\newtheorem{definition}{Definition} \newtheorem{example}{Example} 
\newtheorem{remark}{Remark}  
\newtheorem{corollary}{Corollary}

\def\dim{{\rm \ dim}\,}

\def\cd{{\rm \ cd} }

\def\hut{{\rm \ ht} }
\def\pd{{\rm \ projdim}\  }
\def\depth{{\rm \ depth}\, }
\def\card{{\rm \ card\,} }
\def\rad{{\rm \ rad\,} }
\def\ara{{\rm \ ara\,} }
\def\mod{{\rm \ mod} }
\def\m{{\bf m} }
\def\reg{{\rm \ reg}\, }

\newcommand{\la}{\langle}
\newcommand{\ra}{\rangle}
\def\bv{\hbox{\bf V}} 
  
\def\bn{\hbox{\it I\hskip -2pt N}}

\def\bp{\hbox{\it I\hskip -2pt P}}

\newcommand{\Fcal}{{\mathcal F}}
\newcommand{\Hcal}{{\mathcal H}}
\newcommand{\Lcal}{{\mathcal L}}
\newcommand{\Ical}{{\mathcal I}}
\newcommand{\Jcal}{{\mathcal J}}

\newcommand{\Acal}{{\mathcal A}}

\newcommand{\Rcal}{{\mathcal R}}

\newcommand{\Pcal}{{\mathcal P}}
\newcommand{\Tcal}{{\mathcal T}}
\newcommand{\Vcal}{{\mathcal V}}
\newcommand{\Mcal}{{\mathcal M}}
\newcommand{\Qcal}{{\mathcal Q}}
\newcommand{\Xcal}{{\mathcal X}}
\newcommand{\Dcal}{{\mathcal D}}

\def\demo{\noindent{\bf Proof \ }} 

\begin{document}

%\date{}
%\maketitle
\begin{center}
\uppercase{{\bf  Simplicial ideals, 2-linear ideals and arithmetical rank }}
\end{center}
\advance\baselineskip-3pt
\vspace{2\baselineskip}
\begin{center}
{{\sc Marcel
Morales}\\
{\small Universit\'e de Grenoble I, Institut Fourier, 
UMR 5582, B.P.74,\\
38402 Saint-Martin D'H\`eres Cedex,\\
and IUFM de Lyon, 5 rue Anselme,\\ 69317 Lyon Cedex (FRANCE)}\\
 }
\vspace{\baselineskip}
\end{center}

{\small \sc Abstract.}{\footnote { first version  may 2006, revised version November 2006}} 
{\small In the first part of this paper we study scrollers and linearly joined varieties.
 Scrollers were introduced in \cite{bm4}, linearly joined varieties are an extenson of scroller and were defined in
 \cite{eghp},  and they proved that  scrollers are defined by homogeneous ideals having
 a $2-$linear resolution.
 A particular class of varieties, of important interest in classical Geometry are 
Cohen--Macaulay varieties of minimal degree,
 they were classified by the successive contribution of Del Pezzo \cite{dp}, Bertini \cite{b},  and Xambo \cite{x}. 
 They appear naturally studying 
the fiber cone of of a codimension two toric ideals \cite{gms1}, \cite{gms2}, \cite{bm1}, \cite{h}, \cite{hm}. Let $S$ be a 
polynomial ring and   $\Ical\subset S$ a homogeneous  ideal defining a sequence of 
linearly-joined varieties.

\begin{itemize}

\item We compute the $\depth S/\Ical$, and the cohomological dimension $\cd(\Ical)$.
 \item We prove that under some hypothesis that 
$c(\Vcal)=\depth S/\Ical -1$, where $c(\Vcal)$ is the connectedness dimension of the algebraic set defined by $\Ical$.
\item We characterize sets of generators of  $\Ical$, and give an effective algorithm to find equations, as an application we 
 prove that $\ara (\Ical)=\pd (S/\Ical)$ in the case where $\Vcal$ is a union of linear spaces, in particular
 this applies to any  square free monomial ideal  having a $2-$ linear resolution.
 \item  In the case where $\Vcal$ is a union of linear spaces, the ideal $\Ical$, can be characterized by a tableau,
 which is an extension of a Ferrer (or Young) tableau.
\item We introduce a new class of ideals called  simplicial ideals, 
ideals defining linearly-joined varieties are a particular case of simplicial ideals.
All these results are new, and extend results in \cite{bm4}, \cite{eghp} .
\end{itemize}}

\vskip 1cm 
\section{Introduction}
Throughout this paper we will work with projective schemes $X\subset \bp^r$, but we adopt the algebraic point of view, 
that is we consider a polynomial ring $S$, graded by its standard graduation and reduced homogeneous ideals.
 Our motivation comes from the study 
of the fiber cone $F(I)$ of a codimension two toric ideal, as it was shown in  \cite{gms1}, \cite{gms2}, \cite{bm1}, 
, \cite{h}, \cite{hm},
$F(I)$ appears to be a Cohen--Macaulay reduced ring having minimal multiplicity. 
Projective algebraic sets arithmetically Cohen--Macaulay with minimal degree were  classified by the successive contribution 
of Del Pezzo \cite{dp}, Bertini \cite{b},  and Xambo \cite{x}, and were characterized homologically by 
Eisenbd-Goto \cite{eg}. In the case of irreducible varieties of minimal degree,
 equations defining such varieties are given, but as was pointed by De Concini-Eisenbud-Procesi \cite{cep}: 
"the precise equations satisfied by
 reducible suvarieties of minimal degree remain mysterious", hence they "stop short of giving 
a normal form for the equations of each type. In \cite{bm2}, \cite{bm4}, we tried to answer to
 this question by describing a set of axioms satisfyied by the ideal of any a linear union of scrolls. Some of our results were extended
  in \cite{eghp}.  
  
  In \cite{bm4}, we have extended the notion of varieties of minimal degree to the notion of scrollers,
 where the algebraic set
 is not assumed to be equidimensional, Eisenbud-Green-Hulek-Popescu \cite{eghp} have defined, more generally, 
the notion of linearly joined varieties,
 without assuming that each irreducible component is a scroll, and they prove that scrollers (I use here our definition)
 are exactly  the $2-$ regular (in the sense of Castelnuovo-Mumford)  projective reduced algebraic sets. 
In this paper we continue to investigated about the structure of scrollers, in the first part of the paper 
we extend the characterization of the ideals of scrollers given in \cite{bm4} to the case of linearly joined varieties,
 as a consequence we can compute some invariants of linearly joined varieties, as the depth,
 the connectdness dimension and the arithmetical rank, as a corrolary we can give an effective algorithm to  describe
 equations of lnearly joined varieties, improving previuos results in  \cite{bm4}, \cite{eghp}.
 As an important corollary we get that for the ideal $\Jcal\subset S$ of any $2-$ regular algebraic
 set wich is a union of linear spaces,
 the arithmetical rank  of $\Jcal$ equals the projective dimension $\pd (S/\Jcal)$. In particular this is true for square
 free monomial ideals having a $2-$ linear resolution. 
Note that this results are independent on the characteristic of the field $K$.
 
 In the second part of this paper, we extend the notion of linearly joined varieties to a linear-union of varieties, 
 A reduced ideal $\Jcal\subset S$ defines a  linear-union of varieties, if $\Jcal\subset S$ is 
the intersection of primes ideals  $\Jcal_i=(\Mcal_i,(\Qcal_i))$ for $i=1,...,l$, where $(\Qcal_i)$ is the ideal 
of some sublinear space,  satisfying the property:
$$\Jcal=(\Mcal_1,...,\Mcal_,\bigcap_{j=1}^s(\Qcal_i)).$$
We define a  class  of linear-union of varieties,  defined by   "Simplicial ideals", a Simplicial ideal is 
 a couple $(\Pcal_G,\tilde \bigtriangleup)$ associated to a simplicial complex $\tilde \bigtriangleup$,
 on a set of vertices $G$, and with facets $G_{1},...,G_{s}$, with some properties.
Recently, in her thesis work, my student Ha Minh Lam \cite{h} has studied a class of Simplicial ideals, which variety is an intersection of
 scrolls, and she has proved that they are  scrollers and have a 2-linear resolution. She also has studied the reduction number 
for some class of Simplicial ideals.
We apply and extends some results of the first part of the paper to Simplicial ideals. 
In fact the methods developped here apply to 
a more general setting, this is part of my work in progress.

The author thanks F. Arslan, M. Barile, Ha Minh Lam,  A. Thoma ,
 S. Yassemi, U. Walther, and  Rashid and Rahim Zaare Nahandi for useful discussions. 

\section{Linearly joined varieties}
 \begin{definition} (See \cite{eghp},\cite{bm4}) An ordered sequence $\Vcal_1,...,\Vcal_l\subset \bp^r$ of
 irreducible projective subvarieties  is
 {\it linearly  joined} if for any 
$i=1,...,l-1$ we have :
$$ \Vcal_{i+1}\cap (\Vcal_1\cup ... \cup \Vcal_i)= {\rm span}(\Vcal_{i+1})\cap {\rm span}(\Vcal_1\cup ... \cup \Vcal_i)$$
where  ${\rm span}(\Vcal)$ is the smallest linear subspace of $\bp^r$ containing   $\Vcal$.
\end{definition} 
Linearly joined varieties were defined first in \cite{bm4}, assuming that each variety $\Vcal_i$ is a scroll 
(there were called scrollers),
 then it was extended to the general case in \cite{eghp}.
 
Here we follow the algebraic point of view developed in \cite{bm4}.
 Let $\bv$ a $K$ vector space of dimension $r+1$,  $S=K[\bv]$, the polynomial ring corresponding to the projective space 
 $\bp^r$. For any set $Q\subset \bv$ we will denote by $\la Q \ra \subset \bv$  the  
$K$--vector space generated by $Q$ 
and  by  $(Q)\subset S$ the ideal generated by $Q$. For all $m=1,...,l$, the irreducible variety $\Vcal_m$  is
 defined  by the reduced ideal  $\Jcal_m\subset S$. The linear variety
 $\Lcal_m:={\rm span}(\Vcal_{m})$ is defined by an ideal generated by independent linear forms, so 
let $\Qcal_m\subset \bv$ be the
 linear space such that   $(\Qcal_m)$ is the ideal defining $\Lcal_m$. We can write $\Jcal_m=(\Mcal_m, (\Qcal_m))$
 where $\Mcal_m$ is an ideal.

By \cite[page 163]{bm4}, to show that the sequence of irreducible projective varieties
 $\Vcal_1,...,\Vcal_l\subset \bp^r$ 
 is  linearly joined, ( we will say that the sequence of ideals   $\Jcal_1,...,\Jcal_l$ is 
 {\it linearly joined}), is equivalent to show that for all $k=2,...,l$ :
$$ \leqno{(*) } \hskip 4cm\Jcal_k+\cap_{i=1}^{k-1} \Jcal_i= (\Qcal_k)+(\cap_{i=1}^{k-1} \Qcal_i).$$
It follows from this relation that 
the sequence  $\Lcal_1, ..., \Lcal_l$ is also  {\it linearly joined}.  We denote 
$\Lcal=\Lcal_1\cup ...\cup \Lcal_l$
and $\Qcal:= (\Qcal_1)\cap...\cap(\Qcal_l)$ its difining ideal.
The Theorem 2.1 of \cite[page 163]{bm4}, can be extended to linearly joined varieties:
more precisely 

\begin{definition}let $\Dcal_1=\Qcal_1$, $\Dcal_i:=\bigcap_{j=1 }^i \Qcal_j $.   
For all $i=2,...,r$ let $\la \Delta _{i }\ra$ be a
  linear space such that  $\Dcal_{i-1 }=\Dcal_i\oplus \la \Delta _{i }\ra$, and  let 
$\Pcal_i$ be a linear space such that $\Qcal_i=\Pcal_i\oplus \Dcal_i$.
 
\end{definition}
It follows from the definition that $\Pcal_{1 }=0,\Dcal_l=0, $ $\Dcal_1\supseteq  \Dcal_2\supseteq ...\supseteq \Dcal_l$, is 
a chain of subvector spaces, and $\Dcal_i=\bigoplus_{j=i+1 }^l \la \Delta _{j }\ra$.
 From \cite[page 163]{bm4}, we have that for all $i=2,...,r$ $\Pcal_{i } \cap \Dcal_{i-1 }= 0, $  $\Qcal_i+\Dcal_{i-1 }=
\Dcal_{i } \oplus \la \Delta _{i }\ra\oplus \Pcal_{i }$  

With the notations introduced before we have:
\begin{theorem}{\label {linearly-joined}} The following conditions are equivalent:
\begin{enumerate}
\item  the sequence of ideals   $\Jcal_1,...,\Jcal_l\subset S:=K[\bv]$ is 
 {\it linearly joined}, 
\item For all $i=1,...,l$, there exist   sublinear spaces  $\Dcal_i,\Pcal_i \subset \bv$, with $\Dcal_l=0,\Pcal_1=0$,
 and ideals  $\Mcal_i \subset K[\bv]$ such that 
\begin{itemize}
\item a) for all $i=1,...,l$, $\Jcal_i=(\Mcal_i, \Qcal_i)$ 
\item b) $\Qcal_i =
\Dcal_{i } \oplus  \Pcal_{i }$  
\item c) $\Dcal_1\supseteq  \Dcal_2\supseteq ...\supseteq \Dcal_l$.
\item d) $\Mcal_i\subseteq (\Dcal_{i-1 })$ for all $i=2,...,l$.
\item e) $\Mcal_i\subseteq (\Pcal_{j })$ for all $i=1,...,l-1$ and $j=i+1,...,l$.
\item f) $\bigcap_{j=1 }^{ k-1}(\Qcal_j)\subseteq (\Pcal_k,\Dcal_{k-1 })$ for all $k=2,...,l$.
\end{itemize}
\end{enumerate}
\end{theorem}

\demo Though the proof developped in    \cite[page 164]{bm4} applies here, I  will give  a shorter 
 proof of the implication $" 1. \Rightarrow 2. "$.
 
 The proof  is by induction on $l$. For $l=2$, $\Qcal_{1 }=\la \Delta_{2 }\ra\oplus \Dcal_{2 } , 
\Qcal_{2 }=\Pcal_{2 }\oplus \Dcal_{2 }$ with 
 $\la \Delta_{2 }\ra\cap \Pcal_{2 }=0$, let $S=K[{\bf V}]$, where
 $ {\bf V}\supset \la \Delta_{2 }\ra\oplus \Pcal_{2 }\oplus \Dcal_{2 }$.
 The relation (*) implies that $$(\Mcal_1)+(\Mcal_2)\subset (\la \Delta_{2 }\ra\oplus \Pcal_{2 }\oplus \Dcal_{2 }),$$ and
 without changing the
 ideals $\Jcal_1, \Jcal_2$ we can consider the ideal $\Mcal_1$ modulo the ideal $(\la \Delta_{2 }\ra\oplus \Dcal_{2 }),$
 and the ideal the ideal $\Mcal_2$ modulo the ideal $(\la \Delta_{2 }\ra\oplus \Pcal_{2 }),$ so we get that 
$\Mcal_1\subseteq (\Dcal_{1 }), \Mcal_2\subseteq (\Pcal_{2 })$.

Now suppose that our assertion is true for $l-1$, then 

For $i=1,...,l-1$, $\Qcal_{i }=\Qcal'_{i }\oplus \Dcal_{l-1 }$

$$\Jcal_l+\cap_{i=1}^{l-1} \Jcal_i= (\Qcal_l)+(\cap_{i=1}^{l-1} \Qcal_i)=(\Pcal_{l }\oplus \Dcal_{l-1 }), $$
which implies that $$ \Mcal_1+\Mcal_2+...+\Mcal_l\subset  (\Pcal_{l }\oplus \Dcal_{l-1 }),$$
by induction hypothesis $ \Mcal_1,\Mcal_2,...,\Mcal_{l-1 }$ are defined modulo $(\Dcal_{l-1 })$ and 
without changing the
 ideal $\Jcal_l$ we can consider the ideal $\Mcal_l$ modulo the ideal $ (\Pcal_{l }),$
 and by the same arguments as above we get that for $i=1,...,l-1, \Mcal_i\subseteq (\Pcal_{l }) $, and 
$\Mcal_l\subseteq (\Dcal_{l-1 }) $.
\begin{corollary}

\begin{enumerate}
\item  the sequence of ideals   $\Jcal_1,...,\Jcal_l\subset S:=K[\bv]$ is 
 {\it linearly joined}, 
\item For all $i=1,...,l$, there exist   sublinear spaces  $\Qcal_i\subset \bv$,
 and ideals  $\Mcal_i \subset K[\bv]$ such that 
\begin{itemize}
\item a) for all $i=1,...,l$, $\Jcal_i=(\Mcal_i, \Qcal_i)$ 
\item b) the sequence of ideals   $\Qcal_1,...,\Qcal_l\subset S:=K[\bv]$ is 
 {\it linearly joined}, 
\item c)$\Mcal_i\subseteq (\Qcal_{j })$ for all $i\not=j, i,j\in \{ 1,...,l \}$.
\end{itemize}
\end{enumerate}
\end{corollary}
\demo \begin{enumerate}
\item  "$1.\Rightarrow 2.$", is clear from the above theorem.
\item "$2.\Rightarrow 1.$"
We know from \cite[Prop. 3.4]{eghp}, that $\Jcal_1,...,\Jcal_l\subset S:=K[\bv]$ is 
 {\it linearly joined}, if and only if $\Qcal_1,...,\Qcal_l\subset S:=K[\bv]$ is 
 {\it linearly joined}, and for all $i\not=j, i,j\in \{ 1,...,l \}$ the pair $\Jcal_i,\Jcal_j$ is   
 {\it linearly joined}. So it will be enough to prove that $i\not=j,$  the pair $\Jcal_i,\Jcal_j$ is   
 {\it linearly joined}, but $\Jcal_i+\Jcal_j=(\Qcal_{i }) +(\Qcal_{j })$, by hypothesis so our claim follows.
\end{enumerate}
We also have from of \cite[page 164]{bm4},
\begin{corollary} For any sequence of ideals $\Jcal_1,...,\Jcal_l$ satisfying the axioms a) to e) we have:
$$\bigcap_{j=1 }^{ k}\Jcal_j=(\Mcal_1,...,\Mcal_k, \bigcap_{j=1 }^{ k}(\Qcal_j))$$
for all $k=1,...,l$.

\end{corollary}

\subsection{Equations of linearly-joined Hyperplane arrangements}
As we have seen before if $\Vcal_1, ..., \Vcal_l$ is a sequence of {\it linearly joined}
 irreducible varieties in $\bp^r$, then  
the sequence  $\Lcal_1, ..., \Lcal_l$ is also  {\it linearly joined}. In this subsection we will study this situation.
  We denote 
$\Lcal=\Lcal_1\cup ...\cup \Lcal_l$
and $\Qcal:= (\Qcal_1)\cap...\cap(\Qcal_l)$ its defining ideal.
\begin{corollary}{\label{linearly-joined-cor}}(See Theorem \ref{linearly-joined}) The following conditions are equivalent:
\begin{enumerate}
\item  the sequence of ideals   $(\Qcal_1),...,(\Qcal_l)$ is 
 {\it linearly joined}, 
\item For all $i=1,...,l$, there exist   sublinear spaces  $\Dcal_i,\Pcal_i$, with $\Dcal_l=0,\Pcal_1=0$,
  such that 
\begin{itemize}

\item a) $\Qcal_i =\Dcal_{i } \oplus  \Pcal_{i }$  
\item b) $\Dcal_1\supseteq  \Dcal_2\supseteq ...\supseteq \Dcal_l$.
\item c) $\bigcap_{j=1 }^{ k-1}(\Qcal_j)\subseteq (\Pcal_k,\Dcal_{k-1 })$ for all $k=2,...,l$.
\end{itemize}
\end{enumerate}

\end{corollary}

For all $k=2,...,l$ the sequence  $\Lcal_1, ..., \Lcal_k$ is   {\it linearly joined} in 
the linear space spanned 
by $\Lcal_1, ..., \Lcal_k$. Let
 $\Dcal_{j,k }=\bigoplus_{i=j+1 }^k \la \Delta_{i }\ra $, $\Qcal_{j,k}:=\Pcal_j\oplus \Dcal_{j,k }$.
 So we have that for 
 $j=1,...,k,$   $\Qcal_j=\Qcal_{j,k}\oplus \Dcal_{j,k }$, $\Qcal_{j,k}$ is the ideal defining $\Lcal_j$ in 
 $\la \Lcal_1, ..., \Lcal_k\ra$, and the sequence $\Qcal_{1,k},..., \Qcal_{k,k}$ 
is linearly joined.
 As a consequence of the above corollary we have that 
 $\bigcap_{j=1 }^{ k-1}(\Qcal_{j,k })\subseteq (\Pcal_k)$ for all $k=2,...,l$.
 We can now improve the Lemma 3.1 of \cite{bm4}:
 
 \begin{lemma}\label{equations}For any $k=2,...,l$, $$\bigcap_{j=1 }^{ k}(\Qcal_{j,k })=\bigcup _{j=1 }^{ k} 
(\la \Delta_{j }\ra \times \Pcal_j),$$
where $(\la \Delta_{j }\ra \times \Pcal_j)$ is the ideal generated by all the products $fg$, with 
$f\in\la \Delta_{j }\ra , g\in \Pcal_j.$
\end{lemma}
\demo The proof is by induction on $k$. For $k=2$, $\Qcal_{1,2 }=\la \Delta_{2 }\ra, \Qcal_{1,2 }=\Pcal_{2 }$ and 
 $\la \Delta_{2 }\ra\cap \Pcal_{2 }=0$, and $S=K[\la \Delta_{2 }\ra\oplus \Pcal_{2 }]$. It is clear that 
 $(\la \Delta_{2 }\ra\times \Pcal_{2 })\subset (\Qcal_{1,2 })\cap (\Qcal_{1,2 })$. The other inclusion follows working modulo 
$(\la \Delta_{2 }\ra\times \Pcal_{2 })$ and using the fact that  $\la \Delta_{2 }\ra\cap \Pcal_{2 }=0$.

Now suppose that our assertion is true for $k-1$, then 
$$\bigcap_{j=1 }^{ k}(\Qcal_{j,k })=(\Pcal_k)\cap(\bigcup _{j=1 }^{ k-1} 
(\la \Delta_{j }\ra \times \Pcal_j,\la \Delta_{k }\ra ),$$
but by the above corollary $\bigcup _{j=1 }^{ k-1} 
(\la \Delta_{j }\ra \times \Pcal_j)\subset (\Pcal_k)$, so $\bigcup _{j=1 }^{ k} 
(\la \Delta_{j }\ra \times \Pcal_j)\subset \bigcap_{j=1 }^{ k}(\Qcal_{j,k })$, again using the 
fact that  $\la \Delta_{k }\ra\cap \Pcal_{k }=0$ we will have our statement.
As a consequence of the Lemma we get the following characterization of equations of
 linearly joined hyperplane arrangements:
 
\begin{proposition}{\label{linearly-joined-prop}}(See Theorem \ref{linearly-joined}) The following conditions are equivalent:
\begin{enumerate}
\item  the sequence of ideals   $(\Qcal_1),...,(\Qcal_l)$ is 
 {\it linearly joined}, 
\item For all $i=1,...,l$, there exist   sublinear spaces  $\la \Delta_{i }\ra,\Pcal_i$, with 
$\la \Delta_{1 }\ra=0,\Pcal_1=0$,
  such that 
\begin{itemize}
\item a) $\Dcal_i :=\bigoplus_{j=i+1}^l \oplus \la \Delta_{j }\ra$  
\item b) $\Qcal_i = \Dcal_i\oplus \Pcal_{i }$  
\item c) For any $k=2,...,l$, and $j<k$ we have 
 $$ \la \Delta_{j }\ra \times \Pcal_j \subset (\Pcal_k).$$
\end{itemize}
\demo The above Lemma  implies that for all $k=2,...,l$
 $$\bigcap_{j=1 }^{ k-1}(\Qcal_j)=(\bigcup _{j=2 }^{ k-1} 
(\la \Delta_{j }\ra \times \Pcal_j),\Dcal_{k-1 })\subset (\Pcal_k,\Dcal_{k-1 }).$$
\end{enumerate}

\end{proposition}
\begin{definition}\label{ltilde}
(In view of the applications.) Given a sequence of linearly joined linear spaces $\Lcal_1\cup ...\cup \Lcal_l\subset \bp^r$, 
we will make an extension  $\tilde\Lcal_1\cup ...\cup \tilde\Lcal_l\subset \bp((V\oplus V')^*).$ Let  consider 
a sequence of linear spaces 
$\Hcal_1, ..., \Hcal_l,$ such that $\bp((V\oplus V')^*)=\bp(V^*)\oplus\bigoplus_{i=1}^l  \Hcal_ i$ and set $\tilde\Lcal_i=
\Lcal_i\oplus  \Hcal_ i$.

\end{definition}
\begin{proposition}\label{ltildeprop}\begin{enumerate}
\item The sequence $\tilde\Lcal_1\cup ...\cup \tilde\Lcal_l\subset \bp((V\oplus V')^*)$ is linearly joined.
 Moreover
$$ \tilde\Lcal_i\cap (\tilde\Lcal_1\cup ...\cup \tilde\Lcal_{i-1})=\Lcal_i\cap (\Lcal_1\cup  ...\cup  \Lcal_{i-1}).$$
\item From the algebraic point of view, let $V'=\bigoplus_{i=1}^l  \Fcal_ i$, such that $\Hcal_ i$ is defined by
the ideal $(\bigoplus_{j\not=i}^l  \Fcal_j) $,  so $\tilde\Lcal_i$ is defined by the ideal 
$(\tilde\Qcal_i):=(\Qcal_i\oplus(\bigoplus_{j\not=i}^l  \Fcal_j)$. With this notation
$$\bigcap_{i=1}^l (\tilde\Qcal_i)= (\bigcap_{i=1}^l (\Qcal_i),\bigcup_{i=1}^l \tilde\Qcal_i\times\Fcal_ i). $$

\end{enumerate}
\end{proposition}
\demo\begin{enumerate}
\item For any $i=2,...,l$ we have 
$$ \tilde\Lcal_i\cap (\tilde\Lcal_1\cup ...\cup \tilde\Lcal_{i-1})\subset 
\tilde\Lcal_i\cap \la\tilde\Lcal_1\cup ...\cup \tilde\Lcal_{i-1}\ra $$
 $$= (\Lcal_i\oplus  \Hcal_ i)\cap (\Lcal_1\oplus  \Hcal_ 1\oplus ...\oplus \Lcal_{i-1}\oplus  \Hcal_{i-1})$$
$$\subset \Lcal_i\cap (\Lcal_1+ ...+ \Lcal_{i-1})= \Lcal_i\cap (\Lcal_1\cup  ...\cup  \Lcal_{i-1})$$
The last two relations follows since by hypothesis $\bp((V\oplus V')^*)=\bp(V^*)\oplus\bigoplus_{i=1}^l  \Hcal_ i$, and 
the sequence $\Lcal_1, ..., \Lcal_l\subset \bp(V)$ is linearly joined.
\item  We have the decomposition :
$$ \tilde\Qcal_i=\tilde\Pcal_i\oplus \tilde\Dcal_i,$$ with
$$\tilde\Pcal_i=\Pcal_i\oplus\bigoplus_{j<i}^l  \Fcal_j, \tilde\Dcal_i=\Dcal_i\oplus\bigoplus_{j>i}^l  \Fcal_j.$$
Let $\la  \Delta_{i+1}\ra$ be a   linear space such that $\Dcal_i=\Dcal_{i+1}\oplus\la  \Delta_{i+1}\ra$, so we have that
 $\tilde\Dcal_i=\tilde\Dcal_{i+1}\oplus\la  \Delta_{i+1} \ra\oplus \Fcal_{i+1}, $ and applying the lemma\ref{equations}  the 
 ideal $\bigcap_{i=1}^l(\tilde\Qcal_i)$ is generated by 
 $$\bigcup_{i=2}^l(\la  \Delta_{i} \ra\oplus \Fcal_{i})\times (\Pcal_i\oplus\bigoplus_{j<i}^l  \Fcal_j),$$
 which is equal to:
$$\bigcup_{i=2}^l(\la  \Delta_{i} \ra)\times (\Pcal_i)\cup 
\bigcup_{i=2}^l(\la  \Delta_{i}\ra )\times (\bigoplus_{j<i}^l  \Fcal_j)
\bigcup_{i=2}^l( \Fcal_{i})\times (\Pcal_i)
\bigcup_{i=2}^l( \Fcal_{i})\times (\bigoplus_{j<i}^l  \Fcal_j),$$
but $$\bigcup_{i=2}^l( \Fcal_{i})\times (\bigoplus_{j<i}^l  \Fcal_j)=
\bigcup_{i=2}^l( \Fcal_{i})\times (\bigoplus_{j\not=i}^l  \Fcal_j),$$
and $$\bigcup_{i=2}^l(\la  \Delta_{i}\ra )\times (\bigoplus_{j<i}^l  \Fcal_j)= 
\bigcup_{i=1}^{l-1}(\bigoplus_{j>i}^l\la  \Delta_{i}\ra )\times (  \Fcal_i)=
\bigcup_{i=1}^{l-1}(\Dcal_{i} )\times (  \Fcal_i),$$ so putting both computations together we get our claim.

\end{enumerate}
\subsection{Depth of linearly-joined varieties}
We recall the following important facts
:
\begin{remark} \begin{enumerate}
\item Let $S=K[G]$ be a ring of polynomials over a field $K$, on a set of variables $G$. For any subset $\sigma \subset G$,
with cardinal $\card \sigma $, the local cohomology group $H^{\card \sigma }_\m(K[\sigma ])$, is an $S-$module isomorphic to 
$$K[\sigma ^{-1}]_{\sigma }\simeq 
\displaystyle\bigoplus_{\alpha \in (1,...,1)+\bn^{\card\sigma }}X^{-\alpha }.$$
It then follows that  $H^{\card\sigma }_\m(K[\sigma ])_{-k }=0  $ for $k<\card \sigma $, and 
\hfill\break
$\dim(H^{\card\sigma }_\m(K[\sigma ])_{-\card\sigma })=1.$
\item  Let  $\Jcal\subset S$ be a reduced homogeneous ideal (for the standard grading in the polynomial ring) 
Suppose that  $\Jcal$ doesnot defines a  linear space,  let 
 $h=\depth S/\Jcal$, then  
$H^{h}_m(S/\Jcal)_{-(h-i)} \not=0$ for at least some $i>0$.
\item  For any two ideals $\Jcal_1,\Jcal_2\subset S$ we have the following exact sequence:
$$0\rightarrow S/\Jcal_1 \cap \Jcal_2\rightarrow S/\Jcal_1\oplus  S/\Jcal_2\rightarrow S/(\Jcal_1 + \Jcal_2)\rightarrow 0$$ which gives rise to the 
long exact sequence:
$$\rightarrow H_\m^{h-1}(S/\Jcal_1+\Jcal_2)\rightarrow H_\m^{h}(S/\Jcal_1\cap \Jcal_2)\rightarrow H_\m^{h}(S/\Jcal_1)
\oplus H_\m^{h}(S/\Jcal_2)\rightarrow H_\m^{h}(S/\Jcal_1+\Jcal_2)\rightarrow $$

\end{enumerate}

\end{remark}

 \begin{lemma}{\label {depthsimp}}
Let $\Xcal_1,\Xcal_2\subset \bp^r$ be a {\it linearly  joined} sequence of  projective subchemes 
(having a proper intersection). 
 Let 
$\Jcal_1,\Jcal_2 $ be the (reduced) ideals of definition of $\Xcal_1,\Xcal_2$, $\Lcal_1=span(X_1), \Lcal_2=span(X_2)$
and $(\Qcal_1),(\Qcal_2)$ its defining ideals.   
  Then 
$$\depth S/\Jcal_1 \cap \Jcal_2=\min \{\depth S/\Jcal_1,\depth S/\Jcal_2,\dim S/(\Qcal_1 + \Qcal_2)+1\}.$$
 Moreover since $\dim S/(\Qcal_1 + \Qcal_2)+1\leq \min\{\dim S/\Jcal_1, \dim S/\Jcal_2\}$ we have the particular cases
 \begin{enumerate}\item If $S/\Jcal_2$ is a Cohen--Macaulay ring then 
$$\depth S/\Jcal_1 \cap \Jcal_2=\min \{\depth S/\Jcal_1,\dim S/(\Qcal_1 + \Qcal_2)+1\}.$$
\item If boths $S/\Jcal_1,S/\Jcal_2$ are Cohen--Macaulay rings then 
$$\depth S/\Jcal_1 \cap \Jcal_2=\dim S/(\Qcal_1 + \Qcal_2)+1.$$
\end{enumerate}

  \end{lemma}
\demo 
We have that
$\Jcal_1+\Jcal_2=(Q_1)+(Q_2)$, so $S/\Jcal_1+\Jcal_2$ is isomorphic to a polynomial ring,
 let $h=\dim S/\Jcal_1+\Jcal_2$, and  $q=\min\{\dim S/\Jcal_1, \dim S/\Jcal_2\},$ it follows that 
 $h+1\leq  q$. So the last two assertions follows from the first one.

We have the following exact sequences:
$$ \leqno{(A)}\hskip 2cm 0\rightarrow H_\m^{i}(S/\Jcal_1\cap \Jcal_2)\rightarrow H_\m^{i}(S/\Jcal_1)
\oplus H_\m^{i}(S/\Jcal_2)\rightarrow \rightarrow 0$$
for either $i<h$ or $i>h+1$ and
$$\leqno{(B)}\hskip 2cm \hskip -1cm 0\rightarrow H_\m^{h}(S/\Jcal_1\cap \Jcal_2)\rightarrow H_\m^{h}(S/\Jcal_1)
\oplus H_\m^{h}(S/\Jcal_2)\rightarrow H_\m^{h}(S/\Jcal_1+\Jcal_2)\rightarrow $$
$$\rightarrow H_\m^{h+1}(S/\Jcal_1\cap \Jcal_2)\rightarrow
 H_\m^{h+1}(S/\Jcal_1)
\oplus H_\m^{h+1}(S/\Jcal_2)\rightarrow 0
$$
If $\min \{\depth S/\Jcal_1,\depth S/\Jcal_2\}<h$ or  $\min \{\depth S/\Jcal_1,\depth S/\Jcal_2\}>h$ we have that 
$$\depth S/\Jcal_1\cap \Jcal_2=\min \{\depth S/\Jcal_1,\depth S/\Jcal_2,\dim S/(\Qcal_1+\Qcal_2)+1\},$$ 
 it remains to consider 
the case $\min \{\depth S/\Jcal_1,\depth S/\Jcal_2\}=h$, 
since the above sequence is graded we have for any integer $i>0$ :
$$\hskip -1cm 0\rightarrow H_\m^{h}(S/\Jcal_1\cap \Jcal_2)_{-h+i }\rightarrow H_\m^{h}(S/\Jcal_1)_{-h+i }
\oplus H_\m^{h}(S/\Jcal_2)_{-h+i }\rightarrow H_\m^{h}(S/\Jcal_1+\Jcal_2)_{-h+i }\rightarrow $$
$$\rightarrow H_\m^{h+1}(S/\Jcal_1\cap \Jcal_2)_{-h+i }\rightarrow
 H_\m^{h+1}(S/\Jcal_1)_{-h+i }
\oplus H_\m^{h+1}(S/\Jcal_2)_{-h+i}\rightarrow 0
$$
but since $S/\Jcal_1+\Jcal_2$ is isomorphic to a polynomial ring of dimension $h$,\hfill\break 
 $H_\m^{h}(S/\Jcal_1+\Jcal_2)_{-h+i }=0$, for any integer $i>0$, so we get 
$$  H_\m^{h}(S/\Jcal_1\cap \Jcal_2)_{-h+i }\simeq  H_\m^{h}(S/\Jcal_1)_{-h+i }
\oplus H_\m^{h}(S/\Jcal_2)_{-h+i } $$
Without restriction we can assume for example that $h=\depth S/\Jcal_1$, 
but  \hfill\break $h< \min\{\dim S/\Jcal_1, \dim S/\Jcal_2\}$, this implies that $\Jcal_1$ cannot define a linear space and so
 $H_\m^{h}(S/\Jcal_1)_{-h+i }\not=0$, for at least some $i>0$, wich implies 
$H_\m^{h}(S/\Jcal_1\cap \Jcal_2)_{-h+i }\not=0.$
So $$\depth S/\Jcal_1\cap \Jcal_2 =h=\min \{\depth S/\Jcal_1,\depth S/\Jcal_2,\dim S/(\Qcal_1+\Qcal_2)+1\}.$$ 

\begin{theorem}{\label {regularitydepth}}
Let $\Vcal_1,...,\Vcal_l\subset \bp^r$ be a {\it linearly  joined} sequence of irreducible projective subvarieties. 
 Let 
$\Vcal=\Vcal_1\cup ...\cup \Vcal_l$, $\Jcal$ the (reduced) ideal of definition of $\Vcal$,
 $\Lcal=\Lcal_1\cup ...\cup \Lcal_l$
and $\Qcal:= (Q_1)\cap...\cap(Q_l)$ its defining ideal.
 then 
 \begin{enumerate}
\item $\depth S/\Qcal= \min_{i=1,...,l-1} \{\dim \Lcal_{i+1}\cap (\Lcal_1\cup ... \cup \Lcal_i)\}+2.$
\item $\depth S/\Jcal=\min \{\depth S/\Jcal_1,...,\depth S/\Jcal_l,\depth S/\Qcal\}.$
\item Assume that for all  $i=1,...,l$ the ring $S/\Jcal_i$ is Cohen--Macaulay, then we have $\depth S/\Jcal=
\depth S/\Qcal$.
\end{enumerate}
\end{theorem}
\demo 1.  The proof is by induction on $l$. If $l=2$, both rings $S/(\Qcal_1), S/(\Qcal_2)$ are Cohen--Macaulay, so
 by the above lemma  we have 
$$\depth S/\Qcal_1\cap \Qcal_2 =\dim S/(\Qcal_1+\Qcal_2)+1,$$ and $\dim \Lcal_{1}\cap \Lcal_2=\dim S/(\Qcal_1+\Qcal_2)-1$.

Now suppose that our claim is true for $l-1$ and we will prove it for $l$. We can apply the above lemma to the ideals
 $\bigcap_{i=1}^{l-1}(\Qcal_i),(\Qcal_l)$ so
 $\depth S/\Qcal=\min\{\depth S/\bigcap_{i=1}^{l-1}(\Qcal_i),\dim S/(\Qcal_l)+(\bigcap_{i=1}^{l-1}\Qcal_i)+1  \} ,$
 but  $ \dim \Lcal_{l}\cap (\Lcal_1\cup ... \cup \Lcal_{l-1})+1=\dim S/(\Qcal_l)+(\bigcap_{i=1}^{l-1}\Qcal_i)$ 
and by induction hypothesis 
$$\depth S/\bigcap_{i=1}^{l-1}(\Qcal_i)= \min_{i=1,...,l-2} \{\dim \Lcal_{i+1}\cap (\Lcal_1\cup ... \cup \Lcal_i)\}+2.$$
So the claim follows.

\noindent 2. The proof is by induction on $l$. The case $l= 2$ follows from the above lemma. We suppose that
 our claim is true for $l-1$ and we will prove it for $l$. We can apply the above lemma to the ideals
 $\bigcap_{i=1}^{l-1}(\Jcal_i),(\Jcal_l)$ so
 $$\depth S/\Jcal=\min\{\depth S/\bigcap_{i=1}^{l-1}(\Jcal_i),
\depth S/(\Jcal_l),\dim S/(\Jcal_l)+(\bigcap_{i=1}^{l-1}\Jcal_i)+1  \} ,$$
 but $(\Jcal_l)+(\bigcap_{i=1}^{l-1}\Jcal_i)=(\Qcal_l)+(\bigcap_{i=1}^{l-1}\Qcal_i)$, 
so $ \dim S/(\Jcal_l)+(\bigcap_{i=1}^{l-1}\Jcal_i)=\dim \Lcal_{l}\cap (\Lcal_1\cup ... \cup \Lcal_{l-1})+1,$ 
and by induction hypothesis 
$$\depth S/\bigcap_{i=1}^{l-1}(\Jcal_i)=\min \{\depth S/\Jcal_1,...,
\depth S/\Jcal_l,\depth S/\bigcap_{i=1}^{l-1}(\Qcal_i)\}$$
So by using our claim 1. the claim 2. follows.
 
The proof of the claim 3. follows by the same arguments developed in the proof of the claim 2.

\begin{corollary}\label{depthltilde}
Let  $\Lcal_1\cup ...\cup \Lcal_l\subset \bp^r$ be a sequence of linearly joined linear spaces, and  
consider its extension $\tilde\Lcal_1\cup ...\cup \tilde\Lcal_l\subset \bp((V\oplus V')^*),$ as defined in the Definition \ref{ltilde}. Then 
$$ \depth K[V\oplus V']/\bigcap_{i=1}^{l}(\tilde\Qcal_i)=\depth K[V]/\bigcap_{i=1}^{l}(\Qcal_i).$$
\end{corollary}
From the proof of the Lemma \ref{depthsimp} we also get:
\begin{corollary}\label{cmintersection}
Let $\Jcal_1,...,\Jcal_l \subset S$ be a  sequence  of  ideals.   
  Then \begin{enumerate}\item If  $S/\Jcal_1,S/\Jcal_2,S/(\Jcal_1+\Jcal_2) $ are Cohen--Macaulay rings and 
\hfill\break $\dim S/(\Jcal_1+\Jcal_2)<\min\{ \dim S/\Jcal_1,\dim S/\Jcal_2\}, $ then 
$$\depth S/\Jcal_1 \cap \Jcal_2=\dim S/(\Jcal_1 + \Jcal_2)+1.$$
and $S/\Jcal_1 \cap \Jcal_2$ is Cohen--Macaulay if and only if 
$$\dim S/\Jcal_1 =\dim S/\Jcal_2=\dim S/(\Jcal_1 + \Jcal_2)+1.$$
\item If  $S/\Jcal_1,...,S/\Jcal_l$ are Cohen--Macaulay rings of the same dimension $d$ and for all $i=2,...,l$,
$S/(\Jcal_{i+1} +\bigcap_{j=1}^{i} \Jcal_j)$ is a Cohen--Macaulay ring of dimension $d-1$ then 
$S/(\bigcap_{j=1}^{l} \Jcal_j)$ is a Cohen--Macaulay ring of dimension $d.$ 
\end{enumerate}
\end{corollary}
\begin{remark} The proof of the Lemma \ref{depthsimp}, provides an effective
 way to compute local cohomology modules for linearly joined varieties, it should be interesting 
to study such kind of local cohomology modules.

The next result was proved in \cite{eghp}, we give here a shorter proof.
\begin{theorem}{\label {regularity}}
Let $\Vcal_1,...,\Vcal_l\subset \bp^r$ be a {\it linearly  joined} sequence of irreducible projective subvarieties. 
 Let 
$\Vcal=\Vcal_1\cup ...\cup \Vcal_l$, $\Jcal$ the (reduced) ideal of definition of $\Vcal$,
 $\Lcal=\Lcal_1\cup ...\cup \Lcal_l$
and $\Qcal:= (Q_1)\cap...\cap(Q_l)$ its defining ideal.
 then $$\reg(\Jcal)=\max \{ 2,\reg(\Jcal_1),...,\reg(\Jcal_l)\}.$$
 \end{theorem} 
\demo For $l=2$, our assertion is clear from the following exact sequence  
$$\hskip -1cm 0\rightarrow H_\m^{h}(S/\Jcal_1\cap \Jcal_2)_{-h+j }\rightarrow H_\m^{h}(S/\Jcal_1)_{-h+j }
\oplus H_\m^{h}(S/\Jcal_2)_{-h+j }\rightarrow H_\m^{h}(S/\Jcal_1+\Jcal_2)_{-h+j }\rightarrow $$
$$\rightarrow H_\m^{h+1}(S/\Jcal_1\cap \Jcal_2)_{-h+j }\rightarrow
 H_\m^{h+1}(S/\Jcal_1)_{-h+j }
\oplus H_\m^{h+1}(S/\Jcal_2)_{-h+j }\rightarrow 0
$$
the general case follows by induction on $l$.

\end{remark}

\begin{example} 
Consider  $S=K[a,b,c,x,y,z,u]$ the   ring of polynomials, and 
$$\Jcal_1=(a,b,c); \Jcal_2=(y,z,a,b); \Jcal_3=(x,z-u,b,c); \Jcal_4=(x-u,y-u,a,c).$$
The sequence $\Jcal_1, ...,\Jcal_4 $ is linearly joined:
\begin{itemize}
\item $\Jcal_1+\Jcal_2=(y,z,a,b,c)$
\item $\Jcal_3+\Jcal_1\cap \Jcal_2=(a,b,c,x,z-u)$
\item $\Jcal_1\cap \Jcal_2\cap\Jcal_3=(b,ac,az-au,ax,cz,cy)$, and $\Jcal_4+ \Jcal_1\cap \Jcal_2\cap\Jcal_3=(a,b,c,x-u,y-u)$
\end{itemize}
In this case $\depth(S/\bigcap_{i=1}^{4}\Jcal_i)=3$.

\end{example}
\subsection{Cohomological  dimension}
Let $\Ical\subset S$ be an ideal (in our situation $S$ will be a graded polynomial ring). 
The cohomological dimension $\cd(\Ical)$ is the highest integer $q $ such that $H^q_{\Ical}(S)\not=0$.
In this subsection we compute the cohomological dimension for some  sequences of linearly joined ideals.

 For any two ideals $\Jcal_1,\Jcal_2\subset S$ we have the following exact sequence:
$$0\rightarrow S/\Jcal_1 \cap \Jcal_2\rightarrow S/\Jcal_1\oplus  S/\Jcal_2\rightarrow S/(\Jcal_1 + \Jcal_2)\rightarrow 0$$
 which gives rise to the 
long exact sequence:
$$\rightarrow H_{\Jcal_1 \cap \Jcal_2}^{h-1}(S)\rightarrow H_{\Jcal_1 + \Jcal_2}^{h}(S)\rightarrow H_{\Jcal_1 }^{h}(S)
\oplus H_{\Jcal_2}^{h}(S)\rightarrow H_{\Jcal_1 \cap \Jcal_2}^{h}(S)\rightarrow .$$
\begin{theorem}{\label {cohomologicald}}
Let $\Vcal_1,...,\Vcal_l\subset \bp^r$ be a {\it linearly  joined} sequence of irreducible projective subvarieties. 
 Let 
$\Vcal=\Vcal_1\cup ...\cup \Vcal_l$, $\Jcal_i$ (resp. $\Jcal$) the (reduced) ideal of definition of $\Vcal_i$
 (resp. $\Vcal$),
 $\Lcal=\Lcal_1\cup ...\cup \Lcal_l$
and $\Qcal:= (Q_1)\cap...\cap(Q_l)$ its defining ideal. We assume that each $\Jcal_i$ is a stci.
 \begin{enumerate}
\item $\cd(\Jcal)=\max _{i=2,...,l}\{\dim_K (\Pcal_i+\Dcal_{i-1})-1\}.$
\item $\cd(\Jcal)=\cd(\Qcal)=\pd (S/\Qcal)$.

\end{enumerate}
\end{theorem}
Since each $\Jcal_i$ is a stci, we have that  $\Jcal_i$ can
 be generated up to radical by a regular sequence of $\hut (\Jcal_i)$ elements. Let $q_i:=\hut (\Jcal_i)$ then 
 $\cd(\Jcal_i)=q_i$, also by definition of linearly joined ideals, (we use freely the notations of section 2.)
 for $i=2,...,l$ we have $(\Jcal_i)+(\bigcap_{j=1}^{i-1}\Jcal_j)\subset (\Pcal_i+\Dcal_{i-1})$, wich implies that
 $$\cd((\Jcal_i)=\hut (\Jcal_i)< \hut (\Pcal_i+\Dcal_{i-1})=\dim_K (\Pcal_i+\Dcal_{i-1})=\cd(\Pcal_i+\Dcal_{i-1}).$$
 The same relation also implies 
$\cd((\Jcal_1)=\hut (\Jcal_1)< \cd(\Pcal_2+\Dcal_{1}).$
 The proof is by induction:
 For $i=2$,   
let $h=\cd(\Pcal_2+\Dcal_{1})$ we have 
 $$\rightarrow H_{\Jcal_1 \cap \Jcal_2}^{h-1}(S)\rightarrow H_{\Jcal_1 + \Jcal_2}^{h}(S)\rightarrow H_{\Jcal_1 }^{h}(S)
\oplus H_{\Jcal_2}^{h}(S)\rightarrow H_{\Jcal_1 \cap \Jcal_2}^{h}(S)\rightarrow .$$
which implies that $\cd(\Jcal_1 \cap \Jcal_2)=\cd(\Pcal_2+\Dcal_{1})-1$.
Now by induction suppose that $i\geq 3$ and
$$\cd(\bigcap_{j=1}^{i-1}\Jcal_j)=max _{j=2,...,i-1}\{\dim_K (\Pcal_j+\Dcal_{j-1})-1\}.$$
Let $\Ical :=\bigcap_{j=1}^{i-1}\Jcal_j$, we consider the exact sequence:
$$\rightarrow H_{\Jcal_i \cap \Ical}^{h-1}(S)\rightarrow H_{\Jcal_i + \Ical}^{h}(S)\rightarrow H_{\Jcal_i }^{h}(S)
\oplus H_{\Ical}^{h}(S)\rightarrow H_{\Jcal_i \cap \Ical}^{h}(S)\rightarrow .$$
where $H_{\Jcal_i + \Ical}^{j+1}(S)=0, H_{\Jcal_i }^{j}(S)=0$, for $j\geq h$. This implies that
$$\cd(\bigcap_{j=1}^{i}\Jcal_j)=max \{\cd((\bigcap_{j=1}^{i-1}\Jcal_j), \dim_K (\Pcal_i+\Dcal_{i-1})-1\}=
max _{j=2,...,i}\{\dim_K (\Pcal_j+\Dcal_{j-1})-1\}.$$
his completes the induction. The second assertion follows from our proof.

\subsection{Connectedness dimension}
\begin{definition}
We recall the definition of connectedness dimension for a noetherian topological space $T$:
$$c(T)=\min\{\dim Z: Z\subseteq  T, Z {\rm \ is\ closed\ and\ } T\setminus Z {\rm\ is\ disconnected } \}.$$
Let $K$   be an algebraically closed field. For any ideal radical ideal $\Ical$ in a polynomial ring $S$ over $K$, 
we set $c(S/\Ical)$ for the connectedness dimension of the affine subvariety defined by $\Ical$ in  $Spec(S)$.
\end{definition}
Let remark that if $\Vcal\subset \bp^r$, is a projective variety defined by a homogeneous reduced ideal $\Ical\subset S$,
 then $c(S/\Ical)=c(\Vcal)+1$. I this section we will use the notations and the results of \cite{br-hs}.

\begin{theorem}Let $X$ be a   topological space, $A,B\subset X$ two closed subspaces having no relation of inclusion, then
 \begin{enumerate}
\item  $c(A\cup B)\leq   \dim (A\cap B)$;
\item if $A\cap B$ is irreducible then $c(A\cup B)\leq   c(A)$;
\item If $  B$ is irreducible then $c(A\cup B)\geq   \min\{c(A), \dim (A\cap B) \}$;
\item If $  B,A\cap B$ are irreducible then $c(A\cup B)=   \min\{c(A), \dim (A\cap B) \}$.

\end{enumerate}
\end{theorem}
\demo \begin{enumerate}
\item It follows from the relation $(A\cup B)\setminus (A\cap B)=(A\setminus (A\cap B))\cup (B\setminus (A\cap B))$.
\item Let $Z\subset A$ be a closed set such that $A\setminus Z$ is disconnected, it will be enough to prove that 
$(A\cup B)\setminus Z$ is disconnected.

Assume that $(A\cup B)\setminus Z$ is connected, by hypothesis we have that $A\setminus Z=U_1\cup U_2$,
 with $U_1,U_2$ non empty closed sets in $A\setminus Z$, such that $U_1\cap U_2=\emptyset $.
 We have that $(A\cup B)\setminus Z=U_1\cup U_2\cup B\setminus Z $, and $U_1,U_2, B\setminus Z$ are non empty closed sets in 
$(A\cup B)\setminus Z$. If  $U_1\cap B\setminus Z=\emptyset $ then we can write 
$(A\cup B)\setminus Z=U_1\cup (U_2\cup B\setminus Z) $, which proves that  $(A\cup B)\setminus Z$ is disconnected, 
we get the same conclusion if $U_2\cap B\setminus Z=\emptyset .$ So we have that $U_1\cap B\setminus Z\not=\emptyset $ and 
$U_2\cap B\setminus Z\not=\emptyset $, but we have that
$$(A\cap B)\setminus Z=(A\setminus Z)\cap  ( B\setminus Z)=(U_1\cap B\setminus Z)\cap (U_2\cap B\setminus Z)$$
this is a contradiction since $A\cap B$ is irreducible, showing our claim.
\item Let $Z\subset A\cup B$ be a closed set such that $A\cup B\setminus Z$ is disconnected, it will be enough to
 prove that either $\dim Z\geq \dim (A\cap B)$ or $\dim Z\geq c (A)$.
 By hypothesis we have that $(A\cup B)\setminus Z=U_1\cup U_2$,
 with $U_1,U_2$ non empty closed sets in $(A\cup B)\setminus Z$, such that $U_1\cap U_2=\emptyset $, this implies that 
 $A\setminus (A\cap Z)=A\setminus Z=(U_1\cap A)\cup (U_2\cap A)$, if both $U_1\cap A, U_2\cap A$ are non empty,
 this relation implies that $\dim Z\geq c (A)$ and we get our claim. So we can assume that either $U_1\cap A=\emptyset $ or
 $U_2\cap A=\emptyset $. On the other hand we have that $B\setminus Z=(U_1\cap B)\cup (U_2\cap B)$ and by  assumptions $B$ 
is irreducible so we have either $U_1\cap B=\emptyset $ or
 $U_2\cap B=\emptyset $. So we have $A\setminus Z=U_i$ and $B\setminus Z=U_j$ for $\{i,j\}=\{1,2\}$. 
 These last two conditions imply 
 that $A\cap B\subset Z$, since $A\setminus Z$ and $B\setminus Z$ are disjoint, so we get that $\dim Z\geq \dim (A\cap B)$ and our claim is proved.
 
\item Follows  from the preceding assertions.
\end{enumerate}
\begin{theorem}{\label {depthconn}}
Let $X$ be a   topological space,
 $\Acal_1,...,\Acal_l\subset X$ be a  sequence of closed irreducible subspaces such that for any $i=1,...,l-1$ 
the intersection
 $\Acal_{i+1}\cap (\Acal_1,...,\Acal_{i})$ is  irreducible then:
 $$ c(\bigcup_{i=1}^l \Acal_i)=\min_{i=1,...,l-1} \{\dim \Acal_{i+1}\cap (\Acal_1\cup ... \cup \Acal_i)\}.$$
In particular let $\Vcal_1,...,\Vcal_l\subset \bp^r$ be a {\it linearly  joined} sequence of irreducible projective subvarieties, 
$\Lcal_i=<\Vcal_i>$ for $i=1,...,l$. 
 Let 
$\Vcal=\Vcal_1\cup ...\cup \Vcal_l$, $\Jcal$ the (reduced) ideal of definition of $\Vcal$,
 $\Lcal=\Lcal_1\cup ...\cup \Lcal_l$
and $\Qcal:= (Q_1)\cap...\cap(Q_l)$ its defining ideal 
Then \begin{enumerate}
\item $\displaystyle c(\Vcal)=\min_{i=1,...,l-1} \{\dim \Lcal_{i+1}\cap (\Lcal_1\cup ... \cup \Lcal_i)\}=c(\Lcal).$

\item $\displaystyle c(S/\Jcal)=c(S/\Qcal)=\depth S/\Qcal-1.$
\item Assume that for all $i=1,...,l$, $\Vcal_i$ is arithmetically Cohen--Macaulay  then 
$\displaystyle c(S/\Jcal)=\depth S/\Jcal-1$.
\end{enumerate}
 \end{theorem}

\demo  The proof is by induction on $l$. For $l=2$, since $\Acal_1, \Acal_2$ are irreducible, we have from the above theorem 
 that
 $c(\Acal_1 \cup \Acal_2)=\dim \Acal_1 \cap \Acal_2,$
 so our claim follows. 
 
 Suppose that our claim is true for $l-1$. By the claim 4. of the above theorem,  we have that 
 $c(\Acal_1 \cup...\cup\Acal_l)=\min \{c(\Acal_1 \cup...\cup\Acal_{l-1}),
\Acal_{l}\cap (\Acal_1\cup ... \cup \Acal_{l-1}) \}.$
 Our claim follows by using the induction hypothesis.
 
 The claim 1., 2. and 3. are consequences of the Theorem \ref{regularitydepth}.

\begin{remark}  It should be interesting to study ideals $\Ical$ in a polynomial ring $S$, having the property
$c(S/\Ical)=\depth S/\Ical-1$. Let recall that Hartshorne have introduced
 and studied the varieties connected in codimension one.
 \end{remark}

\subsection {Arithmetical rank of linearly-joined linear spaces}
In this section for the computation  of the arithmetical rank we use the following result of Schmitt and Vogel\cite{s-v}:
\begin{lemma}Let $R$ be a commutative ring, with identity. Let $P$ be a finite subset of elements of $R$. Let $P_0,...,P_r$
be subsets of $P$ such that:
\begin{itemize}
\item (i) $\bigcup_{l=0 }^r P_l=P$;
\item (ii) $P_0$ has exactly one element;
\item (iii) If $p$ and $p''$ are different elements of $P_l$ ($0\leq l\leq r$)
 there is an integer $l'$ with $0\leq l'<l$ and an element $p'\in P_{ l'}$ such that
 $(pp'')^m\in (p')$ for  some positive integer $m$.
\end{itemize}
We set $q_l=\sum_{ p\in P_l}p.$ Let $(P)$ be the ideal generated by $P$, then
$$\rad{(P) }=\rad{(q_0,...,q_r) }.$$

\end{lemma}
\begin{theorem}{\label{stci-lj}} Let $\bv$ be a $K-$ vector space of dimension $r+1$, $S=K[\bv]$, the polynomial ring graded 
by the standard graduation and  $\bp^r $ the projectif space associated to $S$. Let 
$\Lcal=\Lcal_1\cup ...\cup \Lcal_l\subset \bp^r$ be a {\it linearly  joined} sequence of sublinear 
spaces, 
$\Qcal_i\subset \bv$ be a linear space defining  $\Lcal_i$, for $i=1,...,l$, 
and $\Qcal:= (\Qcal_1)\cap...\cap(\Qcal_l)$ the defining ideal of $\Lcal.$ We use the notations and results of section 4.1.
 
There exists an ordered subset $x_1,...,x_n$ of $\bigcup_{ 2\leq i\leq l}\Delta_i \cup \Pcal_i $ 
such that we can arrange the generators of $(\Qcal_1)\cap  (\Qcal_2)\cap ...\cap (\Qcal_{l})=
\bigcup_{ 2\leq i\leq l}\Delta_i \times  \Pcal_i$ into a triangle having 
 $L(l ):=\max_{ 2\leq i\leq l}\{\card (\Pcal_i)+ \card (D_{i-1 })-1 \}$ lines, as follows:

$$  x_1x_{1,0} \eqno{(1)} $$
$$ x_2x_{2,0},\ \ x_1x_{1,1} \eqno{(2)}$$
$$ x_3x_{3,0},\ \ x_2x_{2,1},\ \ x_1x_{1,2} \eqno{(3)}$$
$$ .\ .\ . $$
$$ .\ .\ . $$
$$ x_{j}x_{j,0},\ \ x_{j-1}x_{j-1,1},\ \ .\ .\ . \ \ x_{1}x_{1,j-1} \eqno{(j)}$$
$$ .\ .\ . $$
$$ .\ .\ . $$

 satisfying the properties: 
\begin{enumerate}
\item For any positive integers $j,m$   we have $x_{j,0}=x_{m,0}$, in what follows we set $x_{j,0}=x_{n}$
\item All the products containing $x_{n}$ appear   in the left diagonal of the triangle.
 \item All the products containing $x_{1}$ appear  in the right diagonal of the triangle.
\item For any  $x_{i}$ appearing in the left diagonal,  there is no holes in the right diagonal labelled $i$ consisting of 
$x_{i} x_{i,0},...,x_{i} x_{i,s_i}$ and the elements of the set 
$\{x_{i,0},x_{i,1},...,x_{i,s_i}\}$ are all linearly independent. In what follows we will set $X_{i,j}$ be 
the linear space spanned by $\{x_{i,0},...,x_{i,j}\}.$
%\item there is no monomial labelled  $x_{i} x_{i,l}$ for $i>n-\delta -1$
\item for any $m>i$ and $k$ if there are two products $x_{i}x_{i,k-i},$ 
$x_{m}x_{m,k-m}$  we have :
\begin{itemize}
\item If $x_{m}\in X_{i,s_i}$ then there exist some $s<k$ such that $x_{m}\in X_{i,s}$
 
\item If  $x_{m}\not\in X_{i,s_i}$ then there exist some $s<k$ such that $x_{m,k+1-m}\in X_{i,s}$.

\end{itemize}
\end{enumerate}

\end{theorem}

The proof is given by induction on $l$, the number of irreducible components. 
For $l=2$, take any basis $ P_2$ of $ \Pcal_2$, so we can range the elements in $\Delta_2 \times P_2 $
 in the following triangle of $\card(\Delta_2)+ \card( P_2)-1 $ 
lines:
\begin{center}
\includegraphics[height=2.5 in,width=3in]{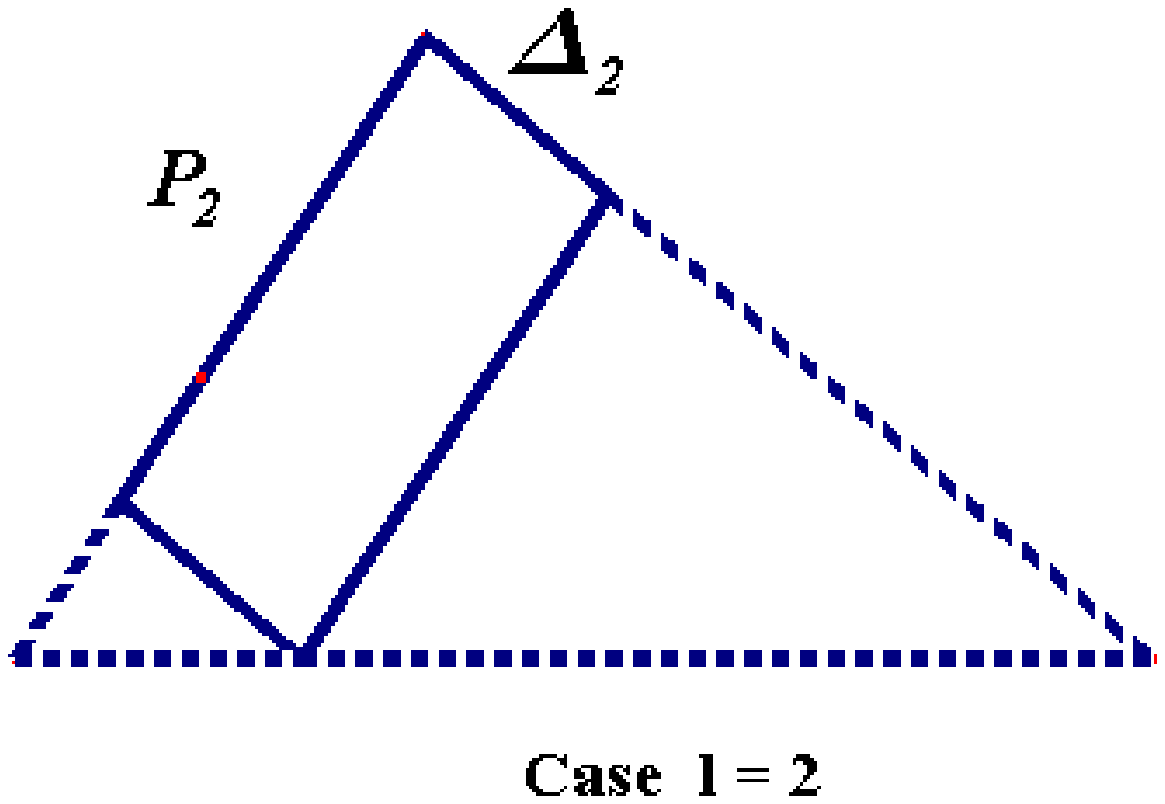}
\end{center}

It is then clear that the theorem is true for $l=2$.

Suppose that the theorem is true for $l-1\geq 2$, and we must prove it for $l$.

The proof is constructive and 
gives an algorithm to find a basis of $ \Pcal_l$ and to compute $\ara((\Qcal_1)\cap  (\Qcal_2)\cap ...\cap (\Qcal_l) ) $.
 
By definition of $D_{l-1} $, we can write $(\Qcal_{i})=(\Qcal'_{i},D_{l-1 })$
 for $i=1,...,l-1$. By induction hypothesis the  generators of $(\Qcal'_1)\cap  (\Qcal'_2)\cap ...\cap (\Qcal'_{l-1}) $  
are ranged in a triangle, satisfying the theorem, then we will define an ordered basis $P_l $ of $ \Pcal_l$, and 
 we form a new triangle by adding 
the quadratic elements   in $\Delta_l \times P_l $ as a diagonal on the left or the right side of this 
triangle.
\begin{center}
\includegraphics[height=2.3 in,width=4.5in]{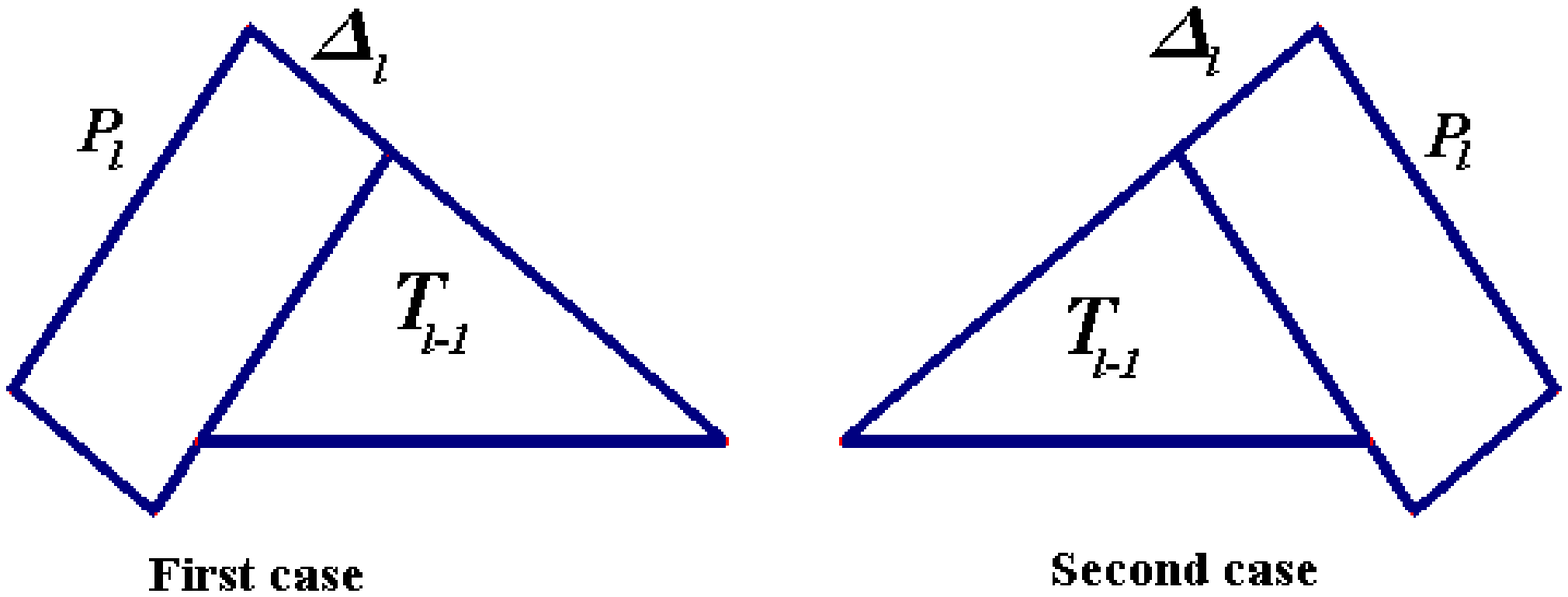}
\end{center}

As a consequence we  have that if  $L(l-1)$ is the number of lines
 in the triangle corresponding to $(\Qcal'_1)\cap  (\Qcal'_2)\cap ...\cap (\Qcal'_{l-1}) $,  
 then:
$$L(l )=\max \{ \card (P_l)+ \card (\Delta  _l)-1, L(l-1 )
+\card (\Delta _l) \} $$ $$= \max_{ 2\leq i\leq l}\{\card (P_i)+ \card (D_{i-1 })-1 \}.$$

 Now let go to the proof. By induction hypothesis there exists an ordered subset $x_1,...,x_n$ of $\bigcup_{ 2\leq i\leq l-1}\Delta_i \cup P_i $ 
such that we can arrange the generators of $(\Qcal'_1)\cap  (\Qcal'_2)\cap ...\cap (\Qcal'_{l-1})=
\bigcup_{ 2\leq i\leq l-1}\Delta_i \times  P_i$ into a triangle of
 $L(l-1)$ lines, as follows:

$$  x_1x_{1,0} \eqno{(1)} $$
$$ x_2x_{2,0},\ \ x_1x_{1,1} \eqno{(2)}$$
$$ x_3x_{3,0},\ \ x_2x_{2,1},\ \ x_1x_{1,2} \eqno{(3)}$$
$$ .\ .\ . $$
$$ .\ .\ . $$
$$ x_{j}x_{j,0},\ \ x_{j-1}x_{j-1,1},\ \ .\ .\ . \ \ x_{1}x_{1,j-1} \eqno{(j)}$$
$$ .\ .\ . $$
$$ .\ .\ . $$

where  $x_{j,0}=x_{n}$  satisfying the properties in the theorem. 

Let $x\in \Delta_l=D_{l-1} $. We consider two cases:
 \begin{itemize}
\item   
$x_{n}\not\in (\Pcal_l)$.     We know that $\bigcup_{ 2\leq i\leq l-1}\Delta_i \times  P_i\subset (\Pcal_l )$,
 in particular for any 
product in the left diagonal $x_i x_n\in (\Pcal_l)$ we have that $x_i \in (P_l)$. Let the set $P_l$ a basis of 
$ \Pcal_l$ 
containing
  the elements appearing in the left diagonal and multiplying $x_n$,
 note that by the point 4 of the theorem they are linearly independent.

We set $x_{n+1}=x$, and we add a left diagonal corresponding to all elements in  $\Delta_l \times  P_l$. 
Then we can range the elements in  $\{x\}\times P_l \cup \bigcup_{ 2\leq i\leq l-1}\Delta_i \times  P_i $
 into the triangle:
$$ x_1 x_{n+1}\eqno{(1)} $$
$$  x_2x_{n+1},\ \ x_1x_n \eqno{(2)} $$
$$ x_3x_{n+1},\ \ x_2x_{n},\ \ x_1x_{1,1} \eqno{(3)}  $$
$$ x_4x_{n+1},\ \ x_3x_{n},\ \ x_2x_{2,1},\ \ x_1x_{1,2} \eqno{(4)}  $$
$$ .\ .\ . $$
$$ .\ .\ . $$
$$ x_jx_{n+1},\ \ x_{j-1}x_{n},\ \ x_{j-2}x_{j-2,1},\ \ .\ .\ .\ \ x_{1}x_{1,j-2} \eqno{(j)}  $$
$$ .\ .\ . $$
$$ .\ .\ . $$
It is clear that we   get the required properties in the theorem.
\item Second  case $ x_n\in (\Pcal_l)$, 
Set $ x_0:=x, x_{0,0}:=x_n$. We define now by induction an ordered basis of  $ \Pcal_l$,
and we add a right diagonal corresponding to all elements in  $\{x\} \times  P_l$.
 Then we can range the elements in  $\{x\}\times P_l \cup \bigcup_{ 2\leq i\leq l-1}\Delta_i \times  P_i $
  in the following triangle:
$$ x_0x_{n}  \eqno{(1)} $$
$$   x_1x_n,\ \ x_0x_{0,1}  \eqno{(2)} $$
$$ x_2x_{n},\ \ x_1x_{1,1},\ \ x_0x_{0,2}  \eqno{(3)}  $$
$$ x_3x_{n},\ \ x_2x_{2,1},\ \ x_1x_{1,2},\ \ x_0x_{0,3}   \eqno{(4)} $$
$$ .\ .\ . $$
$$ .\ .\ . $$
$$ x_{j-1}x_{n},\ \ x_{j-2}x_{j-2,1},\ \ .\ .\ . x_{1}x_{1,j-2},\ \ x_0x_{0,j-1}  \eqno{(j)}  $$
$$ .\ .\ . $$
$$ .\ .\ . $$

Let $\Hcal_s$: Suppose that  we have  defined $x_{0,0},...,x_{0,s}$ lineal independent elements in $ \Pcal_l$,
 and there exists an integer $k>s$, such that for 
any product $x_{l}x_{l,j-l}$ with $l\leq j\leq k$, we have 

 \noindent a) If $x_l\in  \Pcal_l $, then $x_l\in \langle \Pcal_l\rangle_s $, where $\langle \Pcal_l\rangle_s=
\langle x_{0,0},...,x_{0,s}\rangle$;
 
  \noindent  b) If $x_l\not\in  \Pcal_l $ then  $x_{l,j-l}\in \langle \Pcal_l\rangle_s$.
  
We call $k(s)$ the biggest integer $k$ for which  $\Hcal_s$ is true.

The hypothesis $\Hcal_0$ is clearly true.

We suppose that $\Hcal_s$ is  true.
 If $k(s)=L(l-1)+1$ then in order to finish the proof of the proposition, 
we complete to a basis of $ \Pcal_l $.

So suppose that  $k(s)\not=L(l-1)+1$,  let $1\leq i\leq  k(s)$ 
be the smallest integer 
such that  for $x_{i}x_{i,k(s)+1-i}$ the statement in $\Hcal_s$ is not true, so we have two cases:

\begin{itemize}
 \item  If $x_{i}\in  \Pcal_l$,  first we show that necessarily $i=k(s)$. 
Suppose that $i<k(s)$ then $k(s)-i\geq 1$, the element
 $x_{i}x_{i,k(s)-i}$ appears in the line $k(s)$, so by  induction
 hypothesis   $x_i\in \langle \Pcal_l\rangle_s$, in contradiction with the choice of $i$.
 In conclusion $i=k(s)$ and $x_{i}x_{i,k(s)+1-i}=x_{k(s)}x_{k(s),1}$. We set $x_{0,s+1}:=x_{k(s)}$,
 so $\Hcal_{s+1}$ is verified in this
  case
\item If $x_{i}\not\in  \Pcal_l$, then $x_{i,k(s)+1-i}\in  \Pcal_l$, but 
$x_{i,k(s)+1-i}\not\in \langle \Pcal_l\rangle_s$,  we define  $
x_{0,s+1}:=x_{i,k(s)+1-i}$. In order to verify 
$\Hcal_{s+1}$ it will be enough to proof that for any $m>i$, and the element 
$x_{m}x_{m,k(s)+1-m}$  we have either $x_m\in \langle \Pcal_l\rangle_{s+1} $, 
 or $x_m\not\in  \Pcal_l $ and $x_{m,k(s)+1-m}\in \langle \Pcal_l\rangle_{s+1} $. We have to consider several cases

\begin{itemize}\item If $k(s)+1-m=0$, then $x_{k(s)+1,0}=x_{0,0}$, so this case is clear;
\item if $k(s)+1-m\geq 1$, then by induction hypothesis, for the product 
 $x_{m}x_{l,k(s)+1-m}$ we have either $x_m\in X_{i,k(s)}$
 or   $x_{m,k(s)+1-m}\in X_{i,k(s)}$.
 
If  $x_m\in X_{i,k(s)}$, then for each $j\leq k(s)$ we have that $x_{i,j}\in \langle \Pcal_l\rangle_{s}$, it then follows 
that $x_{m}\in \langle \Pcal_l\rangle_{s}$. If  $x_{m,k(s)+1-m}\in X_{i,k(s)}$ 
 we have again that $x_{m,k(s)+1-m}\in \langle \Pcal_l\rangle_{s}$.
The theorem is over.
\end{itemize}
\end{itemize}
\end{itemize}

\begin{theorem} In this theorem $K$ is considered to be algebraically closed.
Let $\Lcal_1,...,\Lcal_l\subset \bp^r$ be a {\it linearly  joined} sequence of linear spaces, 
Let 
$\Qcal_i$ be the  ideal of definition of $\Lcal_i$,
 and $\Qcal:= (\Qcal_1)\cap...\cap(\Qcal_l)$ then :
$$c(S/\Qcal)= \dim S -\ara(\Qcal) -1.$$ 
$$\ara(\Qcal)= \dim S-depth(S/\Qcal).$$
$$\ara(\Qcal)= \pd(S/\Qcal)=\cd(\Qcal).$$
\end{theorem}

Proof. It follows from the Theorem \ref{depthconn} that 
$c(S/\Qcal)=\depth (S/\Qcal)-1$, and  from the above theorem we have that
 $$\ara(\Qcal)\leq \dim S-depth(S/\Qcal)= \dim S-c(S/\Qcal)-1,$$
 on the other hand, by \cite{br-sh}, 19.5.3
 $$c(S/\Qcal)\geq \dim S -\ara(\Qcal) -1,$$
 so we have both  equalities in the claim.
 Recall that by the Auslander-Buchsbaum's theorem $\dim S-depth(S/\Qcal)=\pd(S/\Qcal)$. The last assertion follows from the
  Theorem \ref{cohomologicald}.

\begin{example} Consider again the example 1. $S=K[a,b,c,x,y,z,u]$, and 
$$\Jcal_1=(a,b,c); \Jcal_2=(y,z,a,b); \Jcal_3=(x,z-u,b,c); \Jcal_4=(x-u,y-u,a,c).$$
$\bigcap_{i=1}^{4}\Jcal_i$  is generated by the following terms ordered in a triangle):
$$cb$$
$$ ca   \phantom{cb} \phantom{cb} \phantom{cb} ab$$
$$ cy     \phantom{cb} \hskip 1cm ax  \phantom{cb}  \phantom{cb} b(x-u)$$
$$ cz   \phantom{cb} \phantom{cb}  \phantom{cb}  \phantom{cb}  \hskip 1cm a(z-u)  \phantom{cb}  \phantom{cb} b(y-u)$$
So $\bigcap_{i=1}^{4}\Jcal_i= \rad( cb, ca+ab, cy+ax+b(x-u), cz+a(z-u)+b(y-u)).$

\end{example}

As a corollary we get  the following important theorem.
\begin{theorem} For any square free monomial ideal $\Qcal\subset S$ having a $2-$linear resolution,
 we have $$\ara(\Qcal)= \dim S-depth(S/\Qcal).$$ Moreover computing $depth(S/\Qcal), \ara(\Qcal)$ 
and a set of generators up to radical for $\Qcal$ is effective.

\end{theorem}
Let remark that Herzog-Hibi-Zheng have proved in \cite{hhz} that 
any square free monomial ideal $\Qcal\subset S$ having a $2-$linear resolution, has the property that any power
$\Qcal^k$ has a  linear resolution. In a work in progress we are trying to extend this result to any 
linearly joined hyperplane arrangements.
\subsection{Linearly joined tableau, Ferrer's tableaux.}
\begin{definition}
Suppose that for all $i=1,...,l$, there exist   subsets $ \Delta_{i },P_i$, with 
$ \Delta_{1 }=\emptyset ,P_1=\emptyset $,
  such that $P_i\cap \bigcup_{j=i+1}^l  \Delta_{j }=\emptyset .$
We will say that the elements in $\bigcup_{ 2\leq i\leq l}\Delta_i \times  \Pcal_i$ are ranged in a 
"linearly joined tableau" if
 there exists an ordered subset $x_1,...,x_n$ of $\bigcup_{ 2\leq i\leq l}\Delta_i \cup \Pcal_i $ 
such that we can arrange the elements of  $
\bigcup_{ 2\leq i\leq l}\Delta_i \times  \Pcal_i$ into a triangle  satisfying the properties of the Theorem 
\ref{stci-lj}.
\end{definition}

We have the followig consequence:
\begin{corollary} Any linearly joined sequence of hyperplane arrangements, with ideal $\Qcal\subset S$,
 determines a linearly joined tableau and reciprocally. Moreover 
 $\pd (S/\Qcal)$ is the number of lines in a linearly joined tableau.

\end{corollary}
Note that a linearly joined tableau should be unique up to some operations. This is part of a work in progress.  

The monomial ideals associated to Ferrer's tableaux or diagram are a particular case of the above construction. 
Using the notations of the Corollary 
\ref{linearly-joined-cor}, we can  describe a Ferrer'ideal.

A  Ferrers diagram is a way to represent partitions of a natural number $N$. Let $N,m$ be
 a natural number. A partition of  $N$ is a sum of natural numbers: $N = \lambda _1 + \lambda_2 + ... + \lambda_m$, where 
$\lambda_1 \geq \lambda_2 \geq  ...\geq \lambda_m$. A partition is described by a Young diagram which consists of $m$ rows,
 with the first row containing $\lambda_1$ boxes, 
the second row containing $\lambda_2$ boxes, etc. Each row is left-justified. 
Let $\lambda_{m+1}=0,\delta_0=0,$ and $\delta_1$ be the highest integer such that $\lambda_1=...=\lambda_{\delta_1}, $
 and by induction we define 
$\delta_{i+1}$ as the highest integer such that $\lambda_{\delta_i+1}=...=\lambda_{\delta_{i+1}},$ and set $l$ such that
$\delta_{l-1}=m.$
Let $n= \lambda_1$, we consider two disjoint 
sets of variables :
$\{x_1 , x_2 , ... , x_m\},\{y_1 , y_2 , ... , y_n\}.$ For $i=0,...,l-2$ let
$$\Delta_{l-i}=\{x_{\delta_i+1} ,... , x_{\delta_{i+1}}  \},\Pi_{i+2}=
\{y_{\lambda_{m-i+1}+1}  , ... , y_{\lambda_{m-i} }\}.$$
and $P_{i+2}= \displaystyle\bigcup_{j=2}^{i+2} \Pi_{j}$.  The Ferrer's ideal
 corresponding to  the Ferrer's tableau is generated by 
 $$\Ical_F=(\bigcup_{i=2}^l \Delta_{i}\times P_i).$$
So we have the following:

\begin{proposition}
Let $l>1$ be a natural number  and for $i=1,...,l$ consider two families of  subsets $\Delta_i, \Pi_i$ such that
 $\bv=\bigoplus_{i=2}^l (\la\Delta_i\ra\oplus \la\Pi_i\ra) $ is a decomposition into linear spaces, and let
$$\Pcal_k= \bigoplus_{i=2}^k  \la\Pi_i\ra, \Dcal_{k-1}=\bigoplus_{i=k}^l (\la\Delta_i\ra),$$
$\Qcal_k=\Pcal_k\oplus \Dcal_k,$ and  $(\Qcal_k)\subset K[V]$ be the ideal generated by $\Qcal_k$. Then the linearly joined 
ideal
$\Qcal=\bigoplus_{k=1}^l (\Qcal_k)$ is a Ferrer's ideal.
Reciprocally it is immediate to see that any Ferrer's ideal is obtained in this way. Moreover Ferrer's ideal 
 are characterized as those linearly joined ideals (see Corollary 
\ref{linearly-joined-cor}) $\Qcal=\bigoplus_{k=1}^l (\Qcal_k)$, 
with $\Qcal_k=\Pcal_k\oplus \Dcal_k,$  arrangements of linear spaces 
for which we have the inclusions $\Pcal_{2}\subset \Pcal_{3}\subset ...\Pcal_{l}$ and such that
 $\bv=\Dcal_{1}\oplus \Pcal_{l}$.
\end{proposition}
\demo Applying the algorithm described in the proof of the  Theorem \ref{stci-lj} gives the following tableau, we recognize a Ferrer's tableau, reciprocally any 
Ferrer's tableau gives rise to such decomposition. 
\begin{center}
\includegraphics[height=3 in,width=3.5in]{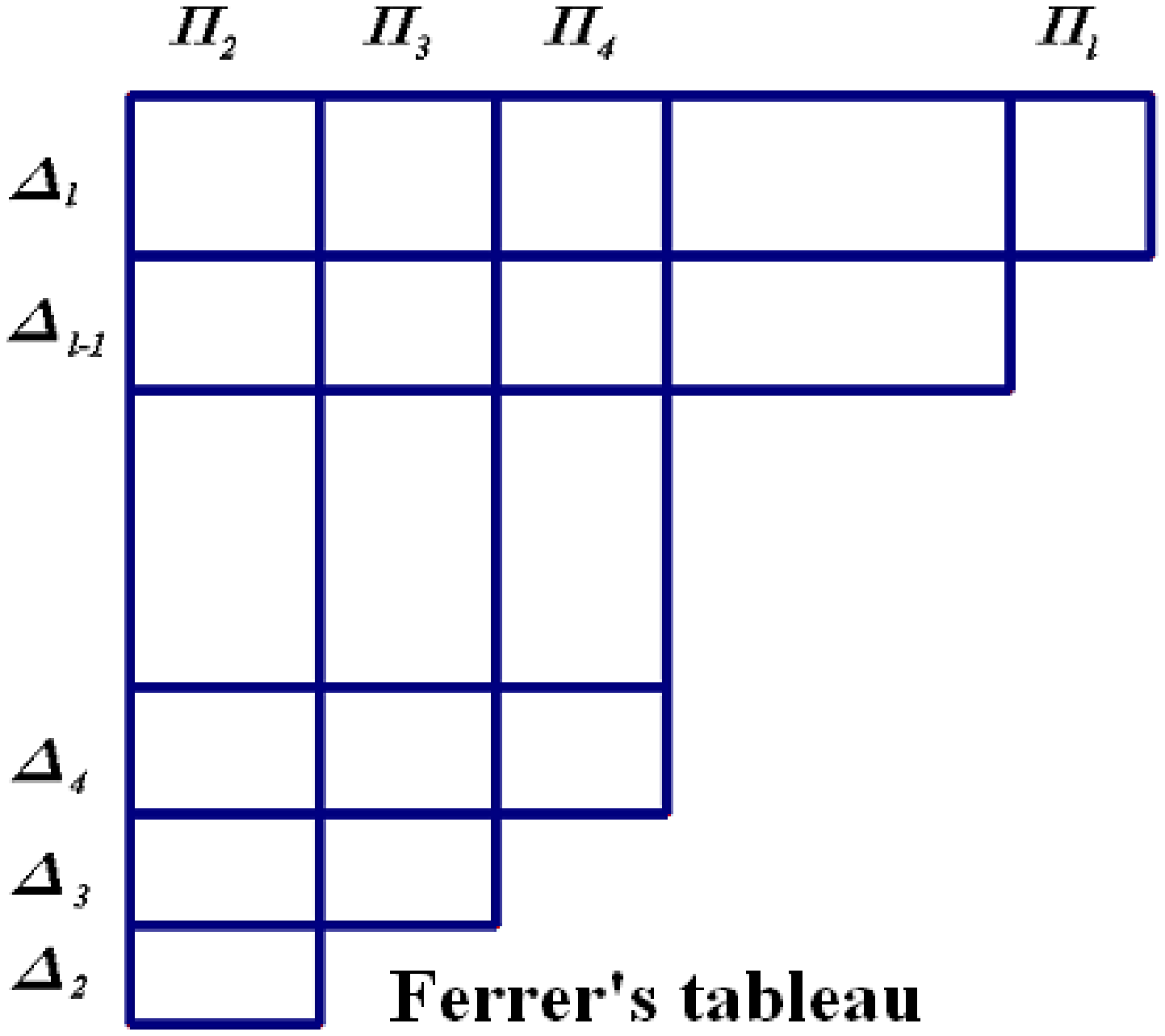}
\end{center}

As a consequence we have 
\begin{corollary}For any Ferrer's ideal $\Ical_\lambda\subset S$,  with label $\lambda_1 \geq \lambda_2 \geq  ...\geq \lambda_m$\begin{enumerate}
\item The minimal primary decomposition of $\Ical_\lambda$ is given inthe above proposition.
 \item $\pd (S/\Ical_\lambda)= \max_{i=1}^m \{\lambda_i+i-1  \}.$
 \item $\ara(\Ical_\lambda)=\cd(\Ical_\lambda)=\pd (S/\Ical_\lambda)$. In fact $\pd (S/\Ical_\lambda)$ is the number of diagonals in a Ferrer's tableau.
\item $c(S/\Ical_\lambda)= \in_{i=1}^m \{(\lambda_m-\lambda_i)+(m-i)  \}.$
\end{enumerate}
 Items 1. and 2. were proved by Corso-Nagel in \cite{cn}.

\end{corollary}

\subsection{Generalized   trees, square free monomial ideals}

Let $K$ be a field, and let $R=K[x_1,...,x_n]$ be the ring of polynomials. Let $\Delta $ be a simplicial complex
 of dimension $d$, on the vertex set $V=\{x_1,...,x_n\}.$ Let $I_\Delta $ be the ideal of $R=K[x_1,...,x_n]$
 generated by the products of those sets of variables which are not faces of $\Delta $. 
The ring $S=K[x_1,...,x_n]/I_\Delta $ is called the Stanley-Reisner ring of $\Delta $ over $K$. It holds that
$\dim S=d+1$. The graph associated to $\Delta $ will be the (1-dimensional) graph $G(\Delta )$ 
on the vertex set $V$ whose edges are the 1-dimensional faces of $\Delta $,
(often named the 1-skeleton of $\Delta $). Vice versa, if $G$ is a graph on the vertex set $V$, 
we shall consider the simplicial complex associated to $G$, denoted by 
$\Delta (G)$, whose maximal faces are all subsets $F$ of $V$ such that the complete graph on $F$ is a subgraph of $G$.
\begin{definition} A generalized $d-$tree on $V$ is a graph  defined recursively as follows:

\begin{itemize}
\item (a) The complete graph on a set of $d+1$ elements of $V$, is a $d-$ tree
\item (b) Let $G$ be a graph on the vertex set $V$. Suppose that there exists some vertex $v\in V$ such that: 
\begin{itemize}
\item the restriction 
$G'$ of $G$ to $V'=V\setminus\{v\}$ is a generalized $d-tree$, and 
\item there exists a subset $V''\subset V'$ of exactly $1\leq j\leq d$ vertices, such that the restriction of $G$ to $V''$
 is a complete graph, and 
\item $G$ is the graph generated by $G'$ and the complete graph on $V''\cup\{v\}$.
\end{itemize}
 
\end{itemize}
The vertex $v$ in the above definition will be called an extremal vertex. If always  $j=d$ in the 
above definition then we say that $G$ is a $d-$tree. In all this paper we will use the terminology generalized 
tree,  instead of generalized $d-$tree.

\end{definition}
We can quote the following theorems of Fröberg \cite{fr}:
\begin{theorem} The Stanley Reisner ring of $\Delta $ is a Cohen-Macaulay ring of minimal degree if and only if
\begin{itemize}
\item the graph $G(\Delta )$ is a $d-$tree and
\item $\Delta =\Delta (G(\Delta )).$
\end{itemize}

\end{theorem}
\begin{theorem} The Stanley Reisner ring of $\Delta $ has a 2-linear resolution if and only if
\begin{itemize}
\item the graph $G(\Delta )$ is a generalized  tree and
\item $\Delta =\Delta (G(\Delta )).$
\end{itemize}

\end{theorem}

Let $\Delta $  a simplicial complex as in the above theorem, 
in the rest of this paper  we will say that $\Delta $  is a generalized tree.
\begin{example}Let $\Ical_1$ be the Stanley-Reisner ideal defined by the simplicial complex:

\begin{center}
\includegraphics[height=2in,width=4in]{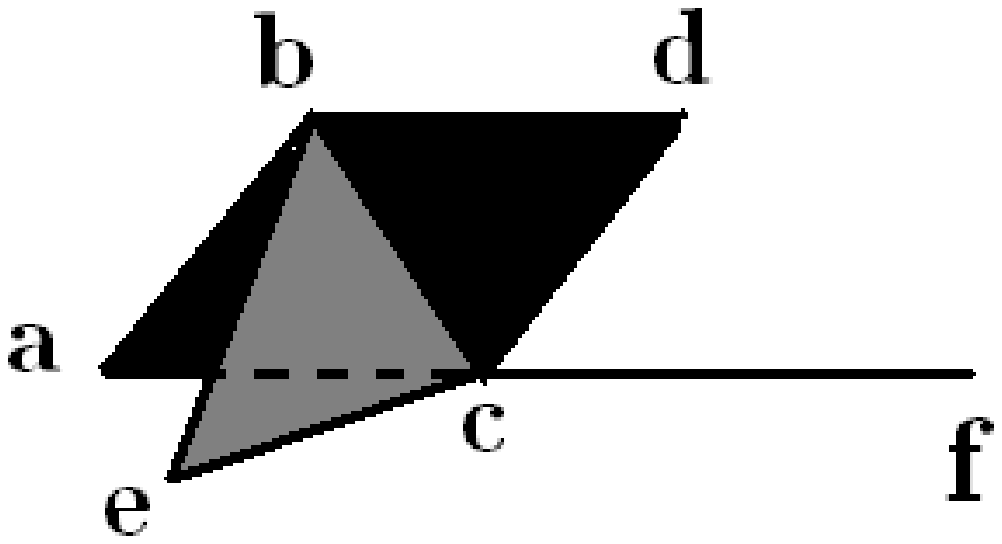}
\end{center}
then the generators of $\Ical_1$  can be ranged in the following linearly joined tableau: 
$$  fd \eqno{(1)} $$
$$ ad,\ \ fe \eqno{(2)}$$
$$ de,\ \ ae,\ \ fb \eqno{(3)}$$
$$ \phantom{ae,\ \ de,\ \ fa ,}\ \ fa \eqno{(4)}$$
\end{example}
\begin{example}Let $\Ical_2$ be the Stanley-Reisner ideal defined by  the simplicial complex:
\begin{center}
\includegraphics[height=2in,width=4in]{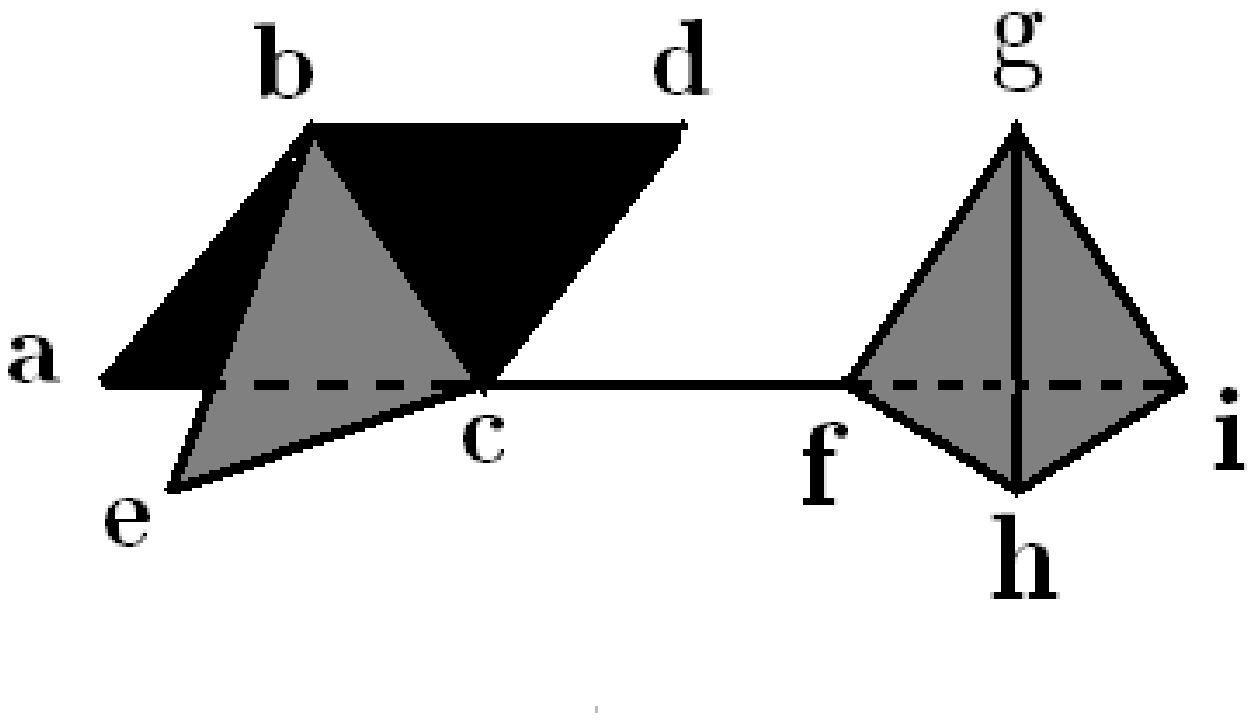}

\end{center}

then the generators of $\Ical_2$ can be ranged in the following linearly joined tableau:
(1st step adding $g$ )
$$  gd \eqno{(0)} $$
$$  fd,\ \ ge \eqno{(1)} $$
$$ ad,\ \ fe,\ \ gb \eqno{(2)}$$
$$ de,\ \ ae,\ \ fb,\ \ ga \eqno{(3)}$$
$$ \phantom{ae,\ \ de,\ \ fa ,}\ \ fa,\ \ gc \eqno{(4)}$$

(2nd step adding $g,h$)
$$  hd \eqno{(-1)} $$
$$  gd,\ \ he  \eqno{(0)} $$
$$  fd,\ \ ge,\ \ hb \eqno{(1)} $$
$$ ad,\ \ fe,\ \ gb,\ \ ha \eqno{(2)}$$
$$ de,\ \ ae,\ \ fb,\ \ ga,\ \ hc \eqno{(3)}$$
$$ \phantom{ae,\ \ de,}\ \ fa,\ \ gc \eqno{(4)}$$
(3th step adding $g,h,i$)
$$  id \eqno{(-2)} $$
$$  hd,\ \ ie \eqno{(-1)} $$
$$  gd,\ \ he,\ \ ib  \eqno{(0)} $$
$$  fd,\ \ ge,\ \ hb,\ \ ia \eqno{(1)} $$
$$ ad,\  fe,\ \ gb,\ \ ha,\ \ ic \phantom{ae,}\eqno{(2)}$$
$$ de,\ \ ae,\ \ fb,\ \ ga,\  hc \phantom{ae,ttttt}\eqno{(3)}$$
$$ \phantom{ae,\ \ de,}\ \ fa,\ \ gc \phantom{ae,ttttt}\eqno{(4)}$$
\end{example}
\begin{example}Let $\Ical_3$ be the Stanley-Reisner ideal defined by  the simplicial complex:
\begin{center}
\includegraphics[height=2in,width=4in]{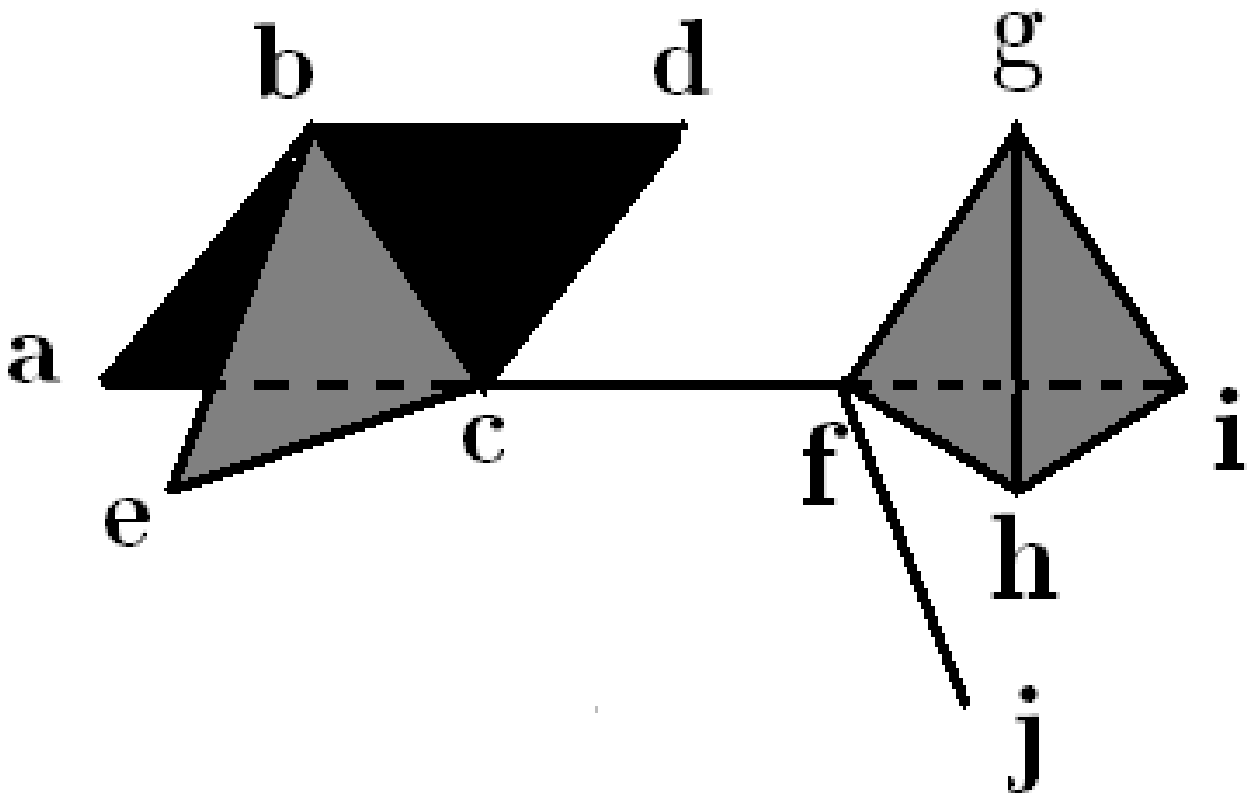}

\end{center}

then the generators of $\Ical_3$ can be ranged in the following linearly joined tableau:

(4th step adding $g,h,i,j$)
$$  jd \eqno{(-2)} $$
$$  id,\ \ ji \eqno{(-2)} $$
$$  hd,\ \ ie,\ \ jh \eqno{(-1)} $$
$$  gd,\ \ he,\ \ ib,\ \ jg  \eqno{(0)} $$
$$  fd,\ \ ge,\ \ hb,\ \ ia,\ \ je \eqno{(1)} $$
$$ ad,\ \ fe,\ \ gb,\ \ ha,\ \ ic,\ \ jb \eqno{(2)}$$
$$ de,\ \ ae,\ \ fb,\ \ ga,\ \ hc,\ \phantom{ae,\ }\ ja \eqno{(3)}$$
$$ \phantom{ae,\ \ de, ttu}\ \ fa,\ \ gc,\ \phantom{ae,\ \ de,}\ jc \eqno{(4)}$$

\end{example}

\section{Linear-union of varieties, simplicial ideals}
\begin{definition}
A reduced  ideal $\Jcal\subset S$ defines a  linear-union of affine varieties, if $\Jcal$ is 
the intersection of primes ideals  $\Jcal_i=(\Mcal_i,(\Qcal_i))$ for $i=1,...,l$, where $(\Qcal_i)$ is the ideal 
of some sublinear space,  satisfying the property:
$$\Jcal=(\Mcal_1,...,\Mcal_,\bigcap_{j=1}^s(\Qcal_i)).$$
\end{definition}

\begin{definition}
Let consider a set of variables $G$ and a decomposition  $G=\displaystyle\bigcup_{i=1}^s  G_{i}$ into distinct  sets  $G_{i}$. 
For any $i=1,...,s$, let $L_i\subset K[G]$ be a set of polynomials, such that $L_i\subset (G_{i})^2$. 
For $i\not=j$, we set  $G_{i,j}=G_{i}\cap G_{j}$, $L_{i,j}=L_i\cap (G_{i,j})$, $I_{i,j}\subset K[G]$ 
be the ideal generated by $L_{i,j}$, $I_{i,i}\subset K[G]$ 
be the ideal generated by $L_{i}\setminus \cup_{j\not= i}L_j$, 
We set  $I_{i}=\displaystyle\sum_{j=1}^s I_{i,j}$,
 $I_G=\displaystyle\sum_{j=1}^s I_{j}$and $\Pcal_G=(I_G,\displaystyle\bigcap_{j=1}^s(G\setminus G_j))$.
We assume 
\begin{itemize}
\item $I_i$ is a prime ideal for any $i$
\item For any $k,l$ and $j \not=k,l$ we have $I_{k,l}\subset (G\setminus G_j).$ 
\end{itemize}We call $\Pcal_G$ a simplicial  ideal.
\end{definition} We can prove our first theorem:
\begin{theorem} Simplicial ideals define linear-union of varieties, more precisely : 
if $\Pcal_G$ is a simplicial  ideal then 
$$\Pcal_G=\bigcap_{j=1}^s(I_j,(G\setminus G_j)).$$

\end{theorem}

We prove the two inclusions:
\begin{itemize}
\item  $"\subset "$: Since 
for $k \not=j,l$ we have that $I_{j,l}\subset (G\setminus G_k)$, it follows that 
$I_G= I_k+\sum_{j,l\not=k} I_{j,l}\subset (I_k+(G\setminus G_k))$,
 on the other hand  for any $k$,
$\cap_{j=1}^s(G\setminus G_j)\subset (G\setminus G_k)$, and so 
$\Pcal_G=(I_G,\cap_{j=1}^s(G\setminus G_j))\subset \cap_{k=1}^s(I_k,(G\setminus G_k)$,
\item $"\supset  "$: Let $\Pcal\supset \Pcal_G$ be a minimal associated prime of
$\Pcal_G$, then $\Pcal\supset \cap_{j=1}^s(G\setminus G_j)$ and since $(G\setminus G_j)$ is a prime ideal,
 there exist some $l$ such that 
$\Pcal\supset (G\setminus G_l)$, on the other hand since $\Pcal\supset \Pcal_G\supset I_l$, it follows that 
$\Pcal\supset (I_l, G\setminus G_l)$, and $(I_l, G\setminus G_l)$ is a prime ideal containing $\Pcal$ by the first item.
 In conclusion
the minimal associated primes of $\Pcal$ are the prime ideals $(I_l, G\setminus G_l)$, for  $l=1,...,s$.

Secondly we compute for  $l=1,...,s$, the $(I_l, G\setminus G_l)$-primary component of  $\Pcal$. In fact we will prove that
$$\Pcal_{(I_l, G\setminus G_l)}= (I_l, G\setminus G_l)_{(I_l, G\setminus G_l)},$$ which will imply that $\Pcal$ is reduced 
and we will get our claim.
Let $j\not= l$, since there exist at least one element  $x\in G_l\setminus G_j$, and $I_{l}\subset (G_{l})^2$, 
we get that $ x\not\in (I_l, G\setminus G_l)$, and $( G\setminus G_j)_{(I_l, G\setminus G_l)}=(1)$, this implies
 that 
 $$\Pcal_{(I_l, G\setminus G_l)}=(I_G,\cap_{j=1}^s(G\setminus G_j))_{(I_l, G\setminus G_l)}=
(I_G,(G\setminus G_l))_{(I_l, G\setminus G_l)}=(I_l,(G\setminus G_l))_{(I_l, G\setminus G_l)},$$
because $I_{k,j}\subset (G\setminus G_l)$ for any $j,k\not=l$, and we are done.
\end{itemize}
\begin{corollary}
$$\dim (\Rcal/\Pcal_G)=\max_{l=1,...,s}\{ \dim (\Rcal/(I_l, G\setminus G_l)) \}=
 \max_{l=1,...,s}\{ \dim (K[G_l]/(I_l))\}.$$
If all our ideals are homogeneous
$$\deg (\Rcal/\Pcal_G)=\sum_{\dim (\Rcal/(I_l, G\setminus G_l))=\dim (\Rcal/\Pcal_G)}
 \deg (\Rcal/(I_l, G\setminus G_l)).$$
\end{corollary}
We illustrate the definition of simplicial ideals by the following examples. In these examples we can apply
 the methods developped above for a linearly joined sequence of ideals, in order to compute $\pd(K[G]/\Pcal_G)$, 
 $\depth(K[G]/\Pcal_G)$, $c(K[G]/\Pcal_G)$ and $cd(\Pcal_G).$
For all of them $K[G]/\Pcal_G$ will be a 
a Cohen--Macaulay ring, by the Corollary\ref{cmintersection}.
We introduce some methods in order to compute the arithmetical rank. We will use it in the next section.
\begin{example} Let $G_1=\{d,b,c,y_1,y_2\}, G_2=\{a,b,c,y_1,y_2, z_1,z_2\}, G_3=\{e,a,c,z_1,z_2\}$ and 
$I_{1,2}$ be the ideal generated by the $2\times 2$ minors of the matrix $M_1$, 
$I_{2,3}$ be the ideal generated by the $2\times 2$ minors of the matrix $M_2$, where 
$$M_1=\pmatrix{b&y_1&y_2\cr y_1&y_2&c\cr } , M_2=\pmatrix{a&z_1&z_2\cr z_1&z_2&c\cr }$$
 we can check easily the hypothesis in the definition of a simplicial ideal. Note 
that $I_{2}:=I_{1,2}+I_{2,3}$ is a prime ideal because it is the toric ideal of the variety parametrized by
$$b=s^3,c=t^3,a=u^3, y_1=s^2t, y_2=st^2, z_1=u^2t, z_2=ut^2,$$ then:
$$\Pcal_G:=(I_{1,2}, I_{2,3}, da,de,be,dz_1dz_2,ey_1,ey_2)=$$
$$(I_{1,2},a,e,z_1,z_2)\cap (I_{1,2},I_{2,3},d,e)\cap (I_{2,3},b,d,y_1,y_2),$$
and $$\Pcal_G=\sqrt{(I_{1,2}, I_{2,3}, da,de,be)},$$
remark that $z_1^2= az_2 \mod I_{2,3}$, so $(dz_1)^2= (da)(dz_2) \mod I_{2,3}$, this shows that
 $dz_1\in \sqrt{(I_{1,2}, I_{2,3}, da,de,be)}$, in the same way we can prove that 
 $dz_2, ey_1,ey_2\in \sqrt{(I_{1,2}, I_{2,3}, da,de,be)}$, which proves the equality proposed.
 
 Let remark that :
 $$\hut (\Pcal_G)=6, ara(I_{1,2})=\hut I_{1,2})=2, ara(I_{2,3})=\hut I_{2,3})=2,$$
 and $( da,de,be)=\sqrt{( de,da+be)}$, it then follows that $\Pcal_G$
 is a stci.
 
\end{example}

\begin{example} Let $G=\{a,b,c,d,e,f,g,h,i,l,m\}$, $G_1=\{d,b,c,f,g,l\}, G_2=\{a,b,c,f,g,h,i\}, G_3=\{e,a,c,h,i,m\}$ and 
$L_{1}$ be the set of   $2\times 2$ minors of the matrix $M_1$, 
$L_{3}$ be the set of   $2\times 2$ minors of the matrix $M_2$, where 
$$M_1=\pmatrix{b&f&g&l\cr f&g&c&d\cr } , M_2=\pmatrix{a&h&i&m\cr h&i&c&e\cr }.$$
%voir fichier macaulay: bin-exemple2
We have that $I_{1}, I_{3}\subset (G_2)$, and we have that $I_{2}:=I_{1,2}+I_{2,3}$ .
 $I_{1}, I_{2},I_{3}$ are prime because they are toric ideals, and
$$\Pcal_G:=(I_{1}, I_{3}, ad,al,be,bm,de,dh,di,dm,ef,eg,el,fm,gm,hl,il,lm)$$ is equal to
$$(I_{1},a,e,m,h,i)\cap (I_{1},I_{3},d,e,l,m)\cap (I_{3},b,d,f,g,l),$$ and
$\pd(K[G]/\Pcal_G)=8.$
On the other hand  $$\Pcal_G=\sqrt{(I_{1}, I_{3}, ad,al,be,bm,de,dm,el,lm)},$$
remark that $h^2= ai \mod I_{3}$, so $(dh)^2= (da)(di) \mod I_{3}$, this shows that
 \hfill\break $dh\in \sqrt{(I_{1}, I_{3}, ad,al,be,bm,de,dm,el,lm)}$, in the same way we can prove our assertion.
 
 Let remark that :
 $$\hut (\Pcal_G)=8, \ara(I_{1})=\hut (I_{1})=3, \ara(I_{3})=\hut (I_{3})=3,$$
 and $( ad,al,be,bm,de,dm,el,lm)=\sqrt{( de,el+md,eb+ml+ad,mb+ld)}$, it then follows
 that $\ara(\Pcal_G)\leq 10.$
 It should be interesting to improve this inequality. Note that $\cd(\Pcal_G)=8, c(K[G]/\Pcal_G)=2.$
\end{example}

\begin{example} Let $G_1=\{d,b,c,f,g,l\}, G_2=\{a,b,c,f,g, h,i\}, G_3=\{e,a,c,h,i,m\}$ and 
$L$ be the set $\{  h^2-ai, f^2-bg, ch-i^2, cf-g^2, bc-fg, ac-hi, b^2d-l^3, ae^2-m^3\}$, this set is
 a generator of the toric 
ideal $I_\Tcal$ parametrized by $u^3-a, s^3-b, t^3-c, v^3-d, w^3-e, s^2t-f, st^2-g, tu^2-h, t^2u-i, s^2v-l, uw^2-m$.
We have that $$L_{1}= \{   f^2-bg,  cf-g^2, bc-fg,  b^2d-l^3\},$$
$$L_{2}= \{  h^2-ai, f^2-bg, ch-i^2, cf-g^2, bc-fg, ac-hi \},$$
$$L_{3}= \{  h^2-ai, ch-i^2,   ac-hi,  ae^2-m^3\}.$$ 
The ideals $I_{1},I_{2},I_{3}$, generated  respectively by $L_{1},L_{2},L_{3}$ are prime, because they are toric,
Then we have that:        
$$\Pcal_G:=(I_{\Tcal}, da,de,dm,dh,di,be,bm,al,ef,eg,el,fm,gm,hl,il,lm)$$ is equal to  
$$(I_{1},a,e,h,i,m)\cap (I_{2},d,e,l,m)\cap (I_{3},b,d,f,g,l),$$ 
On the other hand  $$\Pcal_G=\sqrt{I_{\Tcal}, da,de,be)},$$
remark that $h^2= ah \mod I_{3}$, so $(dh)^2= (da)(di) \mod I_{3}$, this shows that
 $dh\in \sqrt{(I_{\Tcal} , da,de,be)}$, in the same way we can prove that 
 $dm,di,bm,al,ef,eg,el,fm,gm,hl,il,lm\in \sqrt{(I_{\Tcal}, da,de,be)}$, which proves the equality proposed.
Also we have that $\ara(I_{1})\leq 3,\ara(I_{3})\leq 3 $, so $\ara(I_{\Tcal})\leq 6$, which implies that:
 $ 8=\hut(\Pcal_G)\leq \ara(\Pcal_G)\leq \ara(I_{\Tcal})+2\leq 8.$ So $\Pcal_G$ is a
 set theoretically complete intersection and we can give explicitly the generators up to the radical. Note that 
 $I_{\Tcal}$ is also a stci.
 
 Now we study the Cohen--Macaulay property. In this example we have that 
 $$(I_{1},a,h,i,e,m)+ (I_{2},d,l,e,m)=(f^2-bg,  cf-g^2, bc-fg, a,h,i,d,l,e,m) $$ so the quotient ring 
 $K[G]/((I_{1},a,h,i,e,m)+ (I_{2},d,le,m))\simeq K[b,c,f,g]/(f^2-bg,  cf-g^2, bc-fg)$ is Cohen--Macaulay of dimension two,
 this will imply that $K[G]/((I_{1},a,h,i,e,m)\cap  (I_{2},d,le,m))$ is Cohen--Macaulay of dimension three.
 
 $$((I_{1},a,h,i,e,m)\cap  (I_{2},d,l,e,m))+ (I_{3},b,d,f,g,l)=(h^2-ai, ch-i^2,   ac-hi, b,d,e,f,g,l,m) $$
 so the quotient ring 
 $K[G]/(((I_{1},a,h,i,e,m)\cap  (I_{2},d,l,e,m))+ (I_{3},b,d,f,g,l))\simeq K[a,c,h,i]/(h^2-ai, ch-i^2,   ac-hi)$ 
is Cohen--Macaulay of dimension two,
 this will imply from the Corollary \ref{cmintersection} that $K[G]/(\Pcal_G)$ is Cohen--Macaulay of dimension three.
\end{example}

\section{Ara of some simplicial ideals}
In the Definition \ref{ltilde} we have extended any sequence of linearly joined linear spaces. In the case of 
Stanley-Reisner ideal associated to a simplicial complex we can give a more general definition.
  \begin{definition}
Let  $\bigtriangleup(F)$ be a simplicial complex, with set of vertices 
$F$ and facets $F_1,...,F_s$. Let denote by 
$I_{\bigtriangleup(F)},$ the Stanley-Reisner ideal associated to $\bigtriangleup(F)$.
 Let consider a family of disjoints sets
 $F'_l, 1\leq l\leq s $  and disjoints also from $F$, we define  new sets
$G_l= F_l\cup F'_l$, and let $\bigtriangleup(G)$ be the simplicial complex with vertices 
$G=\cup_{1\leq l\leq s }G_l$ and facets $G_1,...,G_s$. We call  $\bigtriangleup(G)$ a extension of $\bigtriangleup(F)$.
\end{definition}
Let remark that  in the above situation $I_{\bigtriangleup(G)},$ 
is generated by $I_{\bigtriangleup(F)},$ and all products $yz$ such that
 $y\in F'_i, y\in G_j\setminus F_i$ for all $i\not=j$.
This is clear since by hypothesis, for all $i\not=j$  and  $y\in F'_i, z\in G_j\setminus F_i$ the edge $[y,z]$ is not in any facet of 
$\bigtriangleup(G)$ so the product $yz$ belongs to $I_{\bigtriangleup(G)}, $ and the other generators of 
$I_{\bigtriangleup(G)} $ have its support in $F$.

It follows from the Proposition \ref{ltildeprop}, that if $\bigtriangleup(F)$ is a generalized tree then 
$\bigtriangleup(G)$ is a generalized tree, and $\depth (K[F]/I_{\bigtriangleup(F)}) =\depth (K[G]/I_{\bigtriangleup(G)})$,

In the following theorem, $I_i$ will be a toric ideal on the variables $G_i$, and we assume that $I_i$ is
 fully parametrized on the set $F_i$, that is for any $y\in F'_i
$, there exists some natural number $m$ such that
 $y^m-x_{F_i}\in I_i$, where $x_{F_i} $ is a monomial with support the set $F_i$. In particular this implies that 
 $\dim K[G_i]/I_i=\card F_i,$ for all $i$. Let remark that the ideal 
$\Pcal_G= (I_1,...,I_s, I_{\bigtriangleup(G)})$ fullfils the conditions to be a simplicial ideal, and it flllows that
$\Pcal_G=\bigcap_{i=1}^s (I_i, G\setminus G_i)$. 
 \begin{theorem} Let consider as above  $\bigtriangleup(H),$  $\bigtriangleup(G)$,  toric ideals 
 $I_i$   on the variables $G_i$, 
  that are fully parametrized on the set $F_i$ and $\Pcal_G=  (I_1,...,I_s, I_{\bigtriangleup(G)})$. 
Then \begin{enumerate}
\item $\Pcal_G=\rad(\sum_{i=1}^s I_i +  I_{\bigtriangleup(F)})$ and $\ara(\Pcal_G)\leq \sum_{i=1}^s \ara(I_i)+ 
\ara I_{\bigtriangleup(F)}$. 
\item If ${\bigtriangleup(F)}$  is 
a generalized  tree and each ideal $I_i$ is a stci then 
$$\cd(\Pcal_G)=\ara(\Pcal_G)= \card G-\depth (K[F]/I_{\bigtriangleup(F)}).$$
$$c(K[G]/\Pcal_G)= \card G-\ara(\Pcal_G) -1.$$ 
\item In particular if ${\bigtriangleup(F)}$  is 
a $d-$tree, and  each ideal $I_i$ is a stci
then $\Pcal_G$ is a stci.
\item If ${\bigtriangleup(F)}$  is 
a generalized  tree,   each ideal $I_i$ is a stci, and $K[G_i]/I_i$ is Cohen--Macaulay, then :
$$c(K[G]/\Pcal_G)= \card G-\ara(\Pcal_G) -1.$$ 
$$\cd(\Pcal_G)=\ara(\Pcal_G)=\pd (K[G]/\Pcal_G).$$
\end{enumerate}

\end{theorem}

Proof. \begin{enumerate}
\item It will be enough to prove that for all $i\not=j$  and  $y\in F'_i, z\in G_j\setminus F_i$, $yz\in \rad(I_1,...,I_s,I_{\bigtriangleup(F)}).$
 Since for every $i$, $I_{i,i}=I_i$ is a simplicial toric ideal, and  $F_i$ is a set of parameters 
of $I_i$, for every element $y\in F'_i
$, there exists some natural number $m$ such that
 $y^m-x_{F_i}\in I_i$, where $x_{F_i} $ is a monomial with support  $F_i$, let $z\in F_j\setminus F_i$,
 it follows that  $zx_{F_i} \in I_{\bigtriangleup(F)}$, which implies that
 $z^my^m\in (I_i,I_{\bigtriangleup(F)})$,  now let 
 let $z\in F'_j\setminus F_i$, then there exists an integer $\mu  $ such that 
 $z^\mu -x_{F_j}\in I_j$, where $x_{F_j} $ is a monomial with support  $F_j\not= F_i$. Let remark that we can choose
$m=\mu $. We have that 
 $$(y^m-x_{F_i})(z^m -x_{F_j})= y^mz^m- y^m x_{F_j}- z^m x_{F_i} + x_{F_i}x_{F_j},$$
 which implies that $y^mz^m\in (I_i,I_j,I_{\bigtriangleup(F)}).$
as a consequence $\Pcal_G=\rad(\sum_i  I_i +  I_{\bigtriangleup(F)})$.
\item If each $I_i$ is a stci then $\ara(I_i)= \card(F'_i)$, and  if $\bigtriangleup(F),$  is 
a generalized  tree, then also $\bigtriangleup(G),$  is 
a generalized  tree, $\depth (K[G]/I_{\bigtriangleup(G)})=\depth (K[F]/I_{\bigtriangleup(F)})$, and
 $\ara(I_{\bigtriangleup(F)})= \card F -\depth (K[F]/I_{\bigtriangleup(F)})$, 
so $$\ara(\Pcal_G)\leq \sum_i \card(F'_i)+ \card F -\depth (K[F]/I_{\bigtriangleup(F)})=
\card G -\depth (K[F]/I_{\bigtriangleup(F)})=$$
$$=\card G -\depth (K[G]/I_{\bigtriangleup(G)})=\card G -c(K[G]/\Pcal_G)-1,$$
since by the Theorem
\ref{depthconn} $\depth (K[G]/I_{\bigtriangleup(G)})=c(K[G]/I_{\bigtriangleup(G)})=c(K[G]/\Pcal_G),$ 
it follows then :
 $$\ara(\Pcal_G)\leq \card G -c(K[G]/\Pcal_G)-1,$$
 or
 $$c(K[G]/\Pcal_G)\leq \card G -\ara(\Pcal_G)-1,$$
 on the other hand, by \cite{br-sh}, 19.5.3
 $$c(K[G]/\Pcal_G)\geq \card G -\ara(\Pcal_G)-1,$$
 so we have the equality.
\item in particular if ${\bigtriangleup(F)}$ is a $d-$tree then $K[F]/I_{\bigtriangleup(F)}$ is a Cohen--Macaulay ring, so
$$\depth (K[F]/I_{\bigtriangleup(F)})=\dim (K[F]/I_{\bigtriangleup(F)})=\dim (K[G]/\Pcal_G)$$
$$\ara(\Pcal_G)\leq \card G -\dim (K[G]/\Pcal_G)=\hut(\Pcal_G),$$ which implies that 
$ \ara(\Pcal_G)=\hut(\Pcal_G), $ and  so $\Pcal_G$ is a stci.
\item  If for every $i$, $K[G_i]/I_i$ is Cohen--Macaulay, then :
$$\depth (K[F]/I_{\bigtriangleup(F)})=\depth (K[G]/I_{\bigtriangleup(G)})=\depth (K[G]/\Pcal_G),
$$ so our statement follows from 2.

\end{enumerate}
The statements 	about cohomological dimension follows from the
  Theorem \ref{cohomologicald}.
\begin{remark}\begin{enumerate}
\item We recall that in \cite{bmt}it was proved that any (simplicial) toric ideal  
 $I$    fully parametrized is a stci if $char(K)=p>0$ and almost stci if $char(K)=p>0$. 
 So we can apply the above theorem to find a large class of examples for
 which $\ara (\Pcal_G)< \pd(K[G]/\Pcal_G).$ This will be published later.
\item It follows from my work \cite{m}, and the work of Robbiano-Valla \cite{rv1},\cite{rv1} that 
 any simplicial  toric ideal 
in codimension two, arithmetically Cohen--Macaulay is a stci. So we can apply 
the above theorem to this family.
\item The examples given in this paper sustend that there is a more general version of the above Theorem,
 this is part of a work in progress.
\item It should be interesting to study homogeneous ideals in a polynomial ring $S$ for  which
 $$c(S/\Ical)= \dim S   -\ara(\Ical)-1.$$
\end{enumerate}
 
\end{remark}

\end{document}